\documentclass[twoside,11pt]{amsart}
\usepackage{amsmath,amssymb}

\newcommand{\CP}{\mathbb{CP}}
\newcommand{\bCP}{\overline{\mathbb{CP}}}
\newcommand{\ta}{\tilde{a}}
\newcommand{\tb}{\tilde{b}}
\newcommand{\tx}{\tilde{x}}
\newcommand{\ty}{\tilde{y}}

\newcommand{\TryPackage}[3]{\IfFileExists{#1.sty}{\usepackage{#1}#2}{#3}
}
\TryPackage{mathrsfs}{\renewcommand{\mathcal}{\mathscr}}{%
        \TryPackage{eucal}{}{}}

\newcommand{\al}{\alpha}
\newcommand{\be}{\beta}

\newcommand{\ZZ}{{\mathbb Z}}
\newcommand{\RR}{{\mathbb R}}

\newcommand{\QQ}{{\mathbb Q}}
\usepackage{latexsym}
\usepackage{amsmath}
\usepackage{amsthm}
\usepackage{amscd}
\usepackage{xypic} %This is XY-pic.
\usepackage{amsfonts}
\usepackage{amssymb}
\usepackage{graphicx}
\usepackage{psfrag}

\newtheorem{df}{Definition}
\newtheorem{thm}[df]{Theorem}

\newtheorem{cor}[df]{Corollary}
\newtheorem{lem}[df]{Lemma}
\newtheorem{prop}[df]{Proposition}

\begin{document}

\title[Constructions of small symplectic 4-manifolds]{Constructions of small symplectic 4-manifolds using Luttinger surgery}

\author{Scott Baldridge}
\author{Paul Kirk}
\date{March 1, 2007}

\thanks{The first   author  gratefully acknowledges support from the
NSF  grant DMS-0507857. The second  author  gratefully acknowledges
support from the NSF  grant DMS-0604310.}

\address{Department of Mathematics, Louisiana State University \newline
\hspace*{.375in} Baton Rouge, LA 70817}
\email{\rm{sbaldrid@math.lsu.edu}}

\address{Mathematics Department, Indiana University \newline
\hspace*{.375in} Bloomington, IN 47405}
\email{\rm{pkirk@indiana.edu}}


\maketitle

\begin{abstract} In this article we use the technique of Luttinger surgery to produce small examples of simply connected and non-simply connected minimal symplectic 4-manifolds.  In particular, we construct:
(1) An example of a minimal symplectic 4-manifold that is homeomorphic but not diffeomorphic to $\CP^2\#3\bCP^2$  which contains a symplectic surface of genus 2, trivial normal bundle, and simply connected complement and a disjoint nullhomologous Lagrangian torus with the fundamental group of the complement generated by one of the loops on the torus.
(2)   A minimal symplectic 4-manifold that is homeomorphic but not diffeomorphic to $3\CP^2\#5\bCP^2$ which has two essential Lagrangian tori with simply connected complement.
 These manifolds can be used to replace $E(1)$ in many known theorems and constructions.  Examples in this article include the smallest known minimal symplectic manifolds with  abelian fundamental groups including  symplectic manifolds with finite and infinite cyclic fundamental group and Euler characteristic 6.
\end{abstract}

%**************************************************
%******* Introduction *****************************
%**************************************************

\section{Introduction}

In this article we construct  a number of small (with respect to the Euler characteristic $e$) simply connected and non-simply connected  symplectic 4-manifolds.

Specifically, we construct examples of:
 \begin{itemize}
\item A minimal symplectic manifold  $X$ homeomorphic but not diffeomorphic to $\CP^2\#3\bCP^2$ containing symplectic genus 2 surface with simply connected complement and trivial normal bundle,  and a disjoint nullhomologous Lagrangian torus (Theorem \ref{cool}).

\item A minimal symplectic manifold  $B$ homeomorphic but not diffeomorphic to $3\CP^2\#5\bCP^2$ containing a disjoint pair of  symplectic tori with  simply connected complement and trivial normal bundle  (Theorem \ref{10baby}).  This provides a smaller substitute for the elliptic surface $E(1)$ in many 4-dimensional constructions.

\item A minimal symplectic manifold $X_1$ with fundamental group $\ZZ$, Euler characteristic $e(X_1)=6$, signature $\sigma(X_1)=-2$ containing a symplectic torus $T$ with trivial normal bundle such that the inclusion $X_1-T\subset X_1$ induces an isomorphism on fundamental groups and so that the inclusion $T\subset X_1$ kills one generator of $\pi_1(T)$ (Theorem \ref{ZZ}). This also provides a smaller substitute for $E(1)$ when only one generator is to be killed.
\end{itemize}

Variations on these constructions quickly provide many more examples of small simply connected minimal symplectic manifolds, including manifolds homeomorphic but not diffeomorphic to  $\CP^2\# 5 \bCP^2$,
$\CP^2\# 7 \bCP^2$,  $3\CP^2\# 7 \bCP^2$, $3\CP^2\# 9 \bCP^2$, and $5\CP^2\# 9 \bCP^2$.  Constructions of small manifolds can also be found in \cite{A,AP,BK2,Kot,park1,PSS,smith, SS,SS1}.

 The manifolds $X,B,X_1$ form building blocks which we use to prove a number of results, including the following.
\begin{itemize}

\item There exists an infinite family of  pairwise non-diffeomorphic smooth simply connected manifolds  each homeomorphic to $\CP^2\#3\bCP^2$.
\item If the group $G$ has a presentation with $g$ generators and $r$ relations, then there exists a symplectic 4-manifold $M$ with fundamental group $G$, $e(M)=10 + 6(g+r)$ and $\sigma(M)=-2(g+r+1)$  (Theorem \ref{50percent}).
\item For any pair of non-negative integers $m,n$ there exists a minimal symplectic manifold homeomorphic but not diffeomorphic to $(1+2m+2n)\CP^2\# (3+6m+4n) \bCP^2$ (Corollary \ref{family}).

\item For any integers $p,q,r$, there exists a symplectic manifold $X_{p,q,r}$ with fundamental group $\ZZ/p\oplus \ZZ/q\oplus \ZZ/r$ with $e=6$ and $\sigma=-2$ (Corollary \ref{abelian}).

\item If the abelian group $G$ is generated by $n$ elements, then there exists a symplectic 4-manifold with fundamental group  $G$, $e = \frac{1}{2} n^2 +\frac{19}{2}n +36 $  and
  $ \sigma =  -\frac{5}{2}n-8.$  (Theorem \ref{genabelian}).

  \item For any non-negative integer $n$, there  exists a symplectic 4-manifold with fundamental group  free of rank $n$, $e = 10 $  and
  $ \sigma = -2.$  (Theorem \ref{free}).

 \item For any symplectic manifold $M$ containing a symplectic surface $G$ of genus 1 or 2  with trivial normal bundle  so that the homomorphism $\pi_1(G)\to \pi_1(M)$ induced by inclusion is trivial, there exists infinitely many smooth manifolds $M_n$ with $e(M_n)=e(M)+ 2+ 4\text{ genus}(G)$, $\sigma(M_n)=\sigma(M)-2$, $\pi_1(M_n)=\pi_1(M)$, and the Seiberg-Witten invariants of $M_n$ are different from those of $M_m$ if $n\ne m$  (Corollary \ref{infinitely2}).

\end{itemize}
 We refer the reader to the body of  the article for more precise statements of these theorems and further results.  One particular feature of our constructions is that they contain  nullhomologous Lagrangian tori for which the method of \cite{FS5} allow us to produce infinitely many non-diffeomorphic but homeomorphic families of manifolds.

Our main tools are Luttinger and torus surgery \cite{Lut, ADK, FS7},  Gompf's symplectic sum construction \cite{Gompf}, and, most importantly,  the Seifert-Van Kampen theorem, which we use to prove our central result, Theorem \ref{prop1}.  This is then combined with Freedman's theorem \cite{Freedman} and fundamental results from Seiberg-Witten theory
\cite{taubes2, taubes3, usher,  FS5, MMS} in the applications.

 A  problem which motivates our investigations concerns uniqueness of the diffeomorphism type of a symplectic manifold which has the smallest Euler characteristic among symplectic manifolds with a fixed fundamental group. For example, for the trivial group, the ``symplectic Poincar\'e conjecture'' (cf. \cite{BK}) asks whether a  symplectic manifold homeomorphic to $\CP^2$ is diffeomorphic to $\CP^2$.  Many constructions on 4-manifolds  are simpler to carry out when the Euler characteristic is large, and this has motivated the problem of finding interesting (e.g.~exotic) simply connected or non-simply connected 4-manifolds with small Euler characteristic. As one works with smaller manifolds, it becomes difficult to alter the smooth structure  without changing the fundamental group or destroying   the existence of a symplectic structure.

Another question which motivates these results concerns whether there is a gap between the Euler characteristics of the best (i.e. smallest) example of a smooth 4-manifold with fundamental group $G$, the best example of a symplectic 4-manifold with fundamental group $G$, and the best example of a complex surface with fundamental group $G$. For example, the smallest   smooth 4-manifold with finite cyclic fundamental group has $e=2$, and there does not exist a smaller smooth manifold.
Corollary \ref{abelian} establishes the existence of a symplectic 4-manifold with finite cyclic fundamental group and $e=6$, this is smallest currently known although it is possible that a smaller one exists. The smallest known complex surface with finite cyclic fundamental group has $e=10$.

The paper is organized as follows.  In Section 2 we describe Luttinger surgery and calculate the fundamental group of the complement of some tori in the 4-torus.  In Section 3 we construct the three main building blocks needed for all subsequent constructions.  In Section 4 we prove our  main result, Theorem \ref{prop1}, which computes the fundamental group (and all meridians and Lagrangian push offs) of  the complement of six Lagrangian tori and a symplectic genus two surface in  a certain symplectic manifold $Z$ satisfying $e(Z)=6,\sigma(Z)=-2$ and $H_1(Z)=\ZZ^6$. With this result in place we construct
 construct the simply connected examples described above and in Section 5 we construct the non-simply connected examples.

To the extent that the methods of the present article focus on quite involved calculations of fundamental groups,  we take great care with our use of the Seifert-Van Kampen theorem,  choice of representative loops, and choices of base points. Some of the fundamental group assertions we prove are perhaps not surprising. However,  the introduction of unwanted conjugation at any stage can easily lead to a loss of control over fundamental groups, in particular leading to plausible but unverifiable calculations.  Given the usefulness of our theorems and that such methods are not so common  in 4-dimensional topology, we feel the care we take is justified.

\bigskip

The authors would like to thank A. Akhmedov, R. Fintushel, C. Judge,  C. Livingston, and J. Yazinsky for helpful discussions.

\section{The fundamental group of the complement of some tori in the 4-torus  }\label{HxK}

\bigskip

\subsection{Luttinger surgery}  Given any Lagrangian torus $T$ in a symplectic 4-manifold $M$, the Darboux-Weinstein theorem
 \cite{MS} implies that there is a parameterization of a tubular neighborhood of $T^2\times D^2\to nbd(T)\subset M$ such that the image of $T^2\times \{d\}$ is Lagrangian for all $d\in D^2$. Choosing any point $d\ne 0$ in $D^2$ gives a push off $F_d:T\to T^2\times \{d\}\subset M-T$ called the {\em Lagrangian push off} or {\em Lagrangian framing}.  Given any embedded curve $\gamma\subset T$, its image $F_d(\gamma)$ is called the {\em Lagrangian push off} of $\gamma$.  The push off depends on the choice of $d$  and the particular parameterization of the tubular neighborhood of $T$,  but its isotopy class in $nbd(T)-T$ depends only on the symplectic structure.   As is common we will abuse terminology slightly and call the isotopy class of $F_d(\gamma)$ for any $d\ne 0$  the Lagrangian push off of $\gamma$.   Any curve isotopic to  $\{t\}\times \partial D^2\subset \partial (nbd(T))$ will be called a {\em meridian} of $T$ and denoted by $\mu$.

  Fix $d\in \partial D$.  If $x,y$ are loops in $T$ generating $H_1(T)$, let $m=F_d(x)$ and $\ell=F_d(y)$.  Then the triple $\mu, m ,\ell$ generate
   $H_1(\partial(nbd(T)))$. Since the 3-torus has abelian fundamental group we may choose a base point $t$ on $\partial(nbd(T))$ and unambiguously refer to $\mu, m, \ell \in\pi_1(\partial(nbd(T)),t)$.

 The push offs and meridians are used to specify coordinates for a {\em $p/q$ torus surgery on $T$ along $\gamma$}. This is   the process of removing a tubular neighborhood of $T$ in $M$ and re-gluing it  so that the embedded curve representing $\mu^p F_d(\gamma)^q$  bounds a disk.  The diffeomorphism type of the resulting manifold depends only on the isotopy class of the identification $T^2\times D^2\to nbd(T)$, and not on the particular point $d$ or the specific choice of $\mu$. Its    fundamental group
is isomorphic to
 \begin{equation}\label{luteq} \pi_1(M-T)/N(\mu^pF_d(\gamma)^q)\end{equation}
 where $N(\mu^pF_d(\gamma)^q)$ denotes the normal subgroup generated by $\mu^pF_d(\gamma)^q$.

When the base point of $M$ is chosen off the boundary of the tubular neighborhood of $T$, the based loops $\mu$ and $\gamma$ are to be joined to the base point {\em by the same path} in $M-T$. Then
Equation (\ref{luteq}) holds with respect to this choice of basing.

Note that if one fixes generating curves $x,y$ on $T$, then the embedded curve $\gamma$ can be expressed in $\pi_1(T)$ in the form $\gamma= x^a y^b$ for some relatively prime pair of integers $ a,b $.  In that case the fundamental group of the manifold obtained by $p/q$ torus surgery on $T$ along $\gamma$
is
$$\pi_1(M-T)/N(\mu^p m^{aq}  \ell^{bq})$$
where, as above, $m=F_d(x)$ and $\ell=F_d(y)$.

The special case of  $p=1, q=k$ is called $1/k$ {\em Luttinger surgery  on $T$ along the embedded curve $\gamma\subset T$}. This yields    a symplectic manifold  (\cite{Lut, ADK}).  The   symplectic form  is unchanged away from a neighborhood of $T$.  The fundamental group of the manifold obtained by $1/k$ Luttinger surgery on $T$ along an embedded curve $\gamma$  is isomorphic to
 $$\pi_1(M-T)/N(\mu F_d(\gamma)^k)$$
 where $N(\mu F_d(\gamma)^k)$ denotes the normal subgroup generated by $\mu F_d(\gamma)^k $.

 It is sometimes convenient to adopt the language of 3-dimensional topology and call the process of gluing $T\times D^2$ to $M-nbd(T)$ a {\em $1/k$ Luttinger filling}, or, more generally, a $p/q$ {\em torus filling}.

 When $p\ne\pm1$ there is no reason why the symplectic form should extend over the neighborhood of $T$, and typically the   smooth manifold obtained by $p/q$ surgery admits no symplectic structure when $p\ne \pm1$.

 \bigskip

\subsection{The complement of two Lagrangian tori in the product of two punctured tori}  Let $\hat H$ and $\hat K$ denote a pair of 2-tori, endowed with the standard symplectic form.    Removing an open disk from $\hat H$ and $\hat K$ yields  punctured tori  $H=\hat H-D$ and $K=\hat K-D$.   View $H\times K$  as a codimension 0 symplectic submanifold  of $T^4=\hat H\times \hat K$ with its  standard product symplectic form. The product $H\times K$ should be considered as the complement of a tubular neighborhood of the (singular) union of two symplectic tori $(\hat H\times \{u_K\})\cup (\{u_H\}\times \hat K)\subset \hat H\times \hat K$ (where $u_H$ and $u_K$ denote the centers of the disks removed.)

Choose a pair of curves $x,y$ representing a standard generating set for $\pi_1(H)$ and a pair of curves $a,b$ representing a standard generating set for $\pi_1(K)$. Let $X,Y$ be parallel push offs of $x$ and $y$ in $H$ and let $A_1,A_2$ be parallel push offs of $a$ in $K$,
as illustrated in the following figure. Let $h$ be the intersection point of $x$ and $y$ and let $k$ be the intersection point of $a$ and $b$.  Give $H\times K$ the base point $(h,k)$.

\bigskip

\begin{figure}[h]
 \begin{center}
\psfrag{y}{$y$} \psfrag{x}{$x$}
\psfrag{X}{$X$} \psfrag{h}{$h$}\psfrag{K}{$K$}\psfrag{H}{$H$}
\psfrag{Y}{$Y$} \psfrag{A1}{$A_1$} \psfrag{A2}{$A_2$} \psfrag{a}{$a$} \psfrag{b}{$b$} \psfrag{k}{$k$}
  \includegraphics[scale=.8]{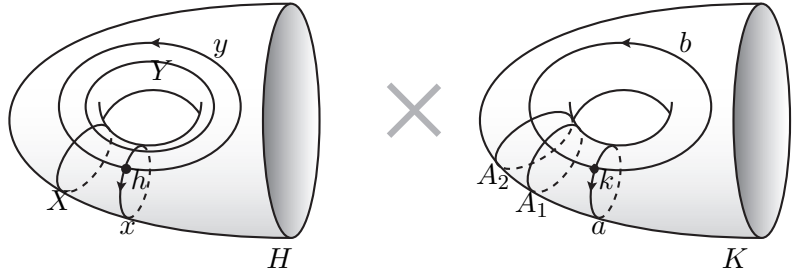}
  \caption{The surface $H\times K$.}
\end{center}
\end{figure}

 \bigskip

We define two disjoint tori $T_1,T_2$ in $H\times K$ as follows.

$$T_1= X\times A_1   \text{ and } T_2= Y \times A_2.$$

  Then the tori $T_1$ and $T_2$ are Lagrangian and the surfaces $H\times \{p\}$ and $
\{q\}\times K$ are symplectic for any $p\in K$ and $q\in H$.

\bigskip

Let     $A_i'$ denote   a  push off of $A_i$  into $K-(A_1\cup A_2)$, $i=1,2$. Then  the parallel tori $T_1'=X\times A_1'$ and $T_2'=Y\times A_2'$  are Lagrangian,
and so the Lagrangian push off  of a curve on $T_i$ is its image in $H\times K-(T_1\cup T_2)$ using this  push off $T_i\to T_i'$.  Sometimes it is preferable to use the push offs using the parallel tori $X'\times A_1$ or $Y'\times A_2$ where $X'$ and $Y'$ are parallel copies in $H$ of $X$ and $Y$.  As we explained above, the manifolds resulting from torus and Luttinger surgery are well defined up to diffeomorphism.

\medskip

The boundary of the tubular neighborhood of $T_i$ is a 3-torus. Therefore $H_1(\partial(nbd(T_i)))=\ZZ^3$, with generating set $\{\mu_i, m_i, \ell_i\}$, where $\mu_i$ is the meridian and  $m_i$ and $\ell_i$ the Lagrangian push offs of  two generators of $H_1(T_i)$.
\bigskip

We specify notation for certain explicit  loops in $H\times K$ based at $(h,k)$.
\begin{enumerate}
\item The loop $x\times\{k\} :I \to  H\times\{k\}$ based at $(h,k)$ will be denoted simply by $x$. This loop misses $T_1\cup T_2$.
\item The loop $y\times\{k\} :I \to  H\times\{k\}$ based at $(h,k)$ will be denoted   by $y$. This loop misses $T_1\cup T_2$.
\item The loop $\{h\}\times a:I\to \{h\}\times K$ based at $(h,k)$ will be denoted by $a$. This loop misses $T_1\cup T_2$.
\item The loop $\{h\}\times b:I\to \{h\}\times K$ based at $(h,k)$ will be denoted by $b$. This loop misses $T_1\cup T_2$.
\end{enumerate}

\bigskip

In \cite[Section 2]{BK2} we proved the following theorem.

\begin{thm} \label{tough} There exist paths in    $H\times K-(T_1\cup T_2)$ from the base point $(h,k)$ to the boundary of the tubular neighborhoods $T_1\times \partial D^2$ and $T_2\times \partial D^2$ with the following property.

Denote  by  $\mu_i, m_i, \ell_i $ the  loops in $H\times K-(T_1\cup T_2)$ based at $  (h,k) $  obtained by following the chosen path to the boundary of the tubular neighborhood of $T_i$, then following (respectively) the meridian of $T_i$ and the two  Lagrangian push offs of the generators  on $T_i$, then returning to the base point along the chosen path.

Then in
$\pi_1(H\times K -(T_1\cup T_2),(h,k))$:
$$\mu_1=[ b^{-1}, y^{-1}], m_1=x, \ell_1= a, $$
and
$$\mu_2=[x^{-1},b], m_2=y, \ell_2=bab^{-1}.$$
where $x,y,a,b$ are  the loops described above.

Moreover,
$\pi_1(H\times K -(T_1\cup T_2),(h,k))$ is generated by $x,y,a,b$ and the relations
$$[x,a]=1, [y,a]=1, [y,bab^{-1}]=1$$
as well as
$$[[x,y],b]=1,[x,[a,b]]=1,[y,[a,b]]=1$$
hold in $\pi_1(H\times K -(T_1\cup T_2),(h,k)).$\qed
\end{thm}

The two  important  things to note in this theorem are, first, the homotopy class of the  loops $x,y,a$, and $b$ based at $(h,k)$ generate
$\pi_1(H\times K -(T_1\cup T_2),(h,k))$. Second, the explicit expressions for $\mu_i,m_i,\ell_i$   allows us to list relations that hold in the fundamental group of the manifold obtained from   torus surgery on the $T_i$ in $H\times K$. For example, the relations
$$[b^{-1},y^{-1}]x^2a^6=1 \text{ and } [x^{-1},b]ba^{-1}b^{-1}=1$$
hold in the fundamental group of manifold the obtained from $H\times K$ by performing $1/2$ surgery on $T_1$ along $m_1\ell_1^3$  and $-1/1$ surgery on $T_2$ along $\ell_2$.

\bigskip
 We will also need the following result.

 \begin{lem} \label{aspherical} Consider the manifold $L$ obtained from $T^4=\hat H\times \hat K$ by performing Luttinger surgeries on $T_1$ along $m_1$ and $T_2$ along either $m_2 $ or $\ell_2$. Then $L$ is aspherical.
 \end{lem}
 \begin{proof} Suppose $L$ is obtained from $T^4$ by performing $1/k_1$ Luttinger surgery on $T_1$ along $m_1$ and  $1/k_2$ Luttinger surgery on $T_2$ along $\gamma$, where $\gamma=m_2$ or $\ell_2$.

 In the  case when $\gamma=m_2$, it is straightforward to see that $L$ is homeomorphic to $Y\times S^1$, where $Y$ is the 3 manifold that fibers over $S^1$ with fiber $\hat H$ and monodromy
 $D_Y^{k_2}\circ D_X^{k_1}$, where $D_X$ and $D_Y$ are the positive Dehn twists along $X$ and $Y$ in $\hat H$. This is explained carefully in \cite[pg. 189]{ADK}.  Thus the universal cover of $L$ is $\RR^4$.

 In the case when $\gamma=\ell_2$, it is not hard to show (see \cite{B1}) that $L$ is homeomorphic to  a non-trivial $S^1$ bundle over $Y$, where $Y$ is the is the 3 manifold that fibers over $S^1$ with fiber $\hat H$ and monodromy
 $D_X^{k_1}$, and the first Chern class of the bundle is $k_2\cdot PD_Y([A_2])$  where $PD_Y$ denotes Poincar\'e duality in $Y$.
 Thus again the universal cover of $L$ is $\RR^4$.

 In either case $L$ is aspherical.
 \end{proof}

 By symmetry, Lemma \ref{aspherical} holds as well if both  surgeries are performed along $\ell_i$.

\section{Three small building blocks}

 \subsection{ }
Our first and simplest building block is the symplectic manifold $W_1=(T^2\times S^2)\# 4\bCP^2$, containing a symplectic genus 2 surface $F_1$ with trivial normal bundle.

We construct the surface $F_1$   by starting with the union of two parallel copies $T^2\times \{p_1\}$, $T^2\times \{p_2\}$ of the torus factor and one copy of $\{q\}\times S^2$ in $T^2\times S^2$. Each of these three surfaces is an embedded symplectic submanifold, and  $\{q\}\times S^2$ intersects each of the tori in one point.
We symplectically resolve the two double points (c.f. \cite{Gompf}), to obtain a symplectic genus 2 surface $F_1$ of square $(2[T]+[S])^2=4$ in $T^4$.    Recall   that topologically, symplectic resolving corresponds to locally replacing a pair of transversely intersecting  discs by an annulus.

Blowing up $T^2\times S^2$ four times at points which lie on $F_1$ and taking the proper transform yields the desired $F_1\subset W_1$. The surface $F_1$ has trivial normal bundle and $W_1$ contains an embedded $-1$ sphere intersecting $F_1$ transversally in exactly one point.

 Let $\phi:F_1\to W_1-nbd(F_1)$ be a push off of $F_1$, and choose a base point $w\in \phi(F_1)$.   Since $F_1$ meets a sphere transversally in one point, the homomorphism induced by inclusion $\pi_1(W_1-F_1,w)\to  \pi_1(W_1,w)$ is an isomorphism.

The two circle coordinates of $T^2$   define classes $s,t\in H_1(W_1)$. Given any base point in $W_1$,  we may unambiguously  write   $\pi_1(W_1)=\ZZ s +\ZZ t$, since $\pi_1(W_1)$ is abelian.

One can choose four loops $s_1,t_1,s_2,t_2$ on $\phi(F_1)$ based at $w$ which generate $\pi_1(\phi(F_1),w)$ and  so that  $[s_1,t_1][s_2,t_2]=1$ in $\pi_1(\phi(F_1),w)$ in such a way that
  the composite $$\pi_1(F_1)\xrightarrow{\phi_*}\pi_1(W_1-nbd(F_1))\cong\pi_1(W_1)\cong  H_1(W_1)$$
takes $s_1$   to $s$, $s_2$ to $s^{-1}$,  $t_1$ to  $t$, and $t_2$ to $t^{-1}$. Thus we adopt the notation:
\begin{enumerate}
\item The loop $s_1:I\to \phi(F_1)\subset \partial (W_1-nbd(F_1))\subset  W_1-nbd(F_1)$ based at $w$ is  a representative loop for the based homotopy class $s\in \pi_1(W_1,w)$.
\item The loop $t_1:I\to \phi(F_1)\subset \partial (W_1-nbd(F_1))\subset  W_1-nbd(F_1)$ based at $w$ is  a representative loop for the based homotopy class $t\in \pi_1(W_1,w)$.
\end{enumerate}

Then the following proposition holds.

\begin{prop} \label{W1} The symplectic surface $F_1$ intersects an embedded sphere transversally in one point and the inclusion $W_1-nbd(F_1)\subset W_1$ induces an isomorphism on fundamental groups. The inclusion $\phi(F_1)\subset W_1-nbd(F_1)$ induces a surjection on fundamental groups.

Moreover,
the loops $s_1, t_1,s_2,t_2$ on $\phi(F_1)$  can be chosen so that
$\pi_1(W_1-nbd(F_1),w)=\ZZ s \oplus \ZZ t$,  where $s,t$ are just the loops $s_1, t_1$ viewed as loops in $W_1-F_1$, and so that the inclusion  $\phi(F_1)\to W_1-nbd(F_1)$ induces the homomorphism $s_1\mapsto s, t_1\mapsto t$, $s_2\mapsto s^{-1}$, and $t_2\mapsto t^{-1}$.  Every $-1$ sphere in $W_1$ intersects $F_1$.
\end{prop}
\begin{proof}
Since   $W_1$ contains a $-1$ sphere intersecting $F_1$ transversally, the meridian of $F_1$ is nullhomotopic in $W_1-nbd(F_1)$.  Hence the inclusion $W_1-nbd(F_1)\to W_1$ induces an isomorphism on fundamental groups by transversality and the Seifert-Van Kampen theorem.  The other assertions about the fundamental group are explained above.

The four exceptional spheres all meet  $F_1$ since the blowup was performed on $F_1$. Denote by $T, S$, $ E_1$, $E_2,E_3,$ and $ E_4$ the five generators of $H_2(W_1)$, where $T=T^2\times \{p\}$, $S=\{q\}\times S$, and the $E_i$ are the exceptional classes. Thus $F_1=2T+S-E_1-E_2-E_3-E_4$. The Hurewicz theorem shows  that the spherical classes are spanned by $S, E_1,E_2,E_3,E_4$.  Consideration of the intersection form shows that  a $-1$ sphere must have the form $aS\pm E_i$.  Then $(aS\pm E_1)\cdot F_1=2a\pm 1\ne 0$. Thus every $-1$ sphere intersects $F_1$.
\end{proof}

Suppose that $P$ is any symplectic 4-manifold containing a symplectic surface $G$ of genus 2 with trivial normal bundle. Then the {\em symplectic sum, $S$,  of $W_1$ and $P$ along $F_1$ and $G$} (c.f. \cite{Gompf}) is a symplectic manifold described topologically as the union of
 $ W_1-nbd(F_1)$ and $ P-nbd(G)$   along their boundary using a fiber preserving diffeomorphism
 $F_1\times S^1\to G\times S^1$ of the boundary of their tubular neighborhoods.    The diffeomorphism type of the manifold $S$ may depend on the choice of such  a diffeomorphism,  which can be specified up to isotopy by choosing trivializations  of the tubular neighborhoods of $F_1$ and $G$ and a diffeomorphism
 $\phi:F_1\to G$. One then glues  $ W_1-nbd(F_1)$ to $ P-nbd(G)$ using the gluing diffeomorphism
 $$\tilde\phi:\partial(W_1-nbd(F_1))=F_1\times S^1\cong \partial(P-nbd(G))
 =G\times S^1, \ (f,z)\mapsto (\phi(f),z).$$
 Then
 $$S=W_1-nbd(F_1)\cup_{\tilde\phi} P-nbd(G).$$
 The symplectic sum is defined more generally when $G$ and $F$ have normal bundles with opposite Euler class, i.e. if $[G]^2=-[F]^2$. For our purposes it will suffice to consider symplectic sums along square zero surfaces. Moreover, the framings we use will either be explicit, or unimportant to the fundamental group calculations.

 \bigskip

 Assume that the base point $p$ of $P$ lies on $G$, and that
 $\phi: F_2\to G$ is base point preserving, $\phi(w)=p$.
Denote by $N$ the subgroup   of $\pi_1(P,p)$ normally generated by
 $\phi(s_1s_2), \phi(t_1t_2),$ and $\phi([s_1,t_1])$. Then Proposition \ref{W1} and the Seifert-Van Kampen theorem imply that
  $\pi_1(S,w) $ is a quotient of $ \pi_1(M,p)/N.$

 \bigskip

\noindent{\bf Remark.}  To properly understand this assertion, it is important to remember that the meridian of $G$ in $P-nbd(G)$  bounds a disk in $S$, namely the punctured exceptional sphere in $W_1-nbd(F_1)$. Thus one can think of attaching
$W_1-nbd(F_1)$ to $ P-nbd(G)$ in two steps: first attach a 2-disk (the punctured exceptional sphere), and then attach the rest of $W_1-nbd(F_1)$.

The first step recovers
$\pi_1(P,p')$ from $\pi_1(P -nbd(G),p')$, where $p'=(p,1)=\tilde\phi (w)\in G\times S^1=\partial(nbd(G))$ is the push off of $p$. Since the homomorphism $\pi_1(F_1\times \{1\}, w)\to \pi_1(W_1-  F_1,w)$ induced by a push off is surjective, the Seifert-Van Kampen theorem together with this observation shows that
$$\pi_1 (P-nbd(G),p')\to \pi_1(S,w)$$
is a surjection which factors through $\pi_1(P,p')$.  Thus one obtains a surjection $\pi_1(P,p')\to \pi_1(S,w)$ whose kernel contains the homotopy classes of the loops  $\tilde \phi(s_1s_2), \tilde \phi(t_1t_2),$ and $\tilde \phi([s_1,t_1])$.   A small path in the normal disk fiber to $G$ identifies $\pi_1(P,p')$ with $\pi_1(P,p)$ and the push off $G\mapsto G\times \{1\}\subset G\times S^1$ is isotopic  to the inclusion $G\subset P$ (using the radial coordinate in $D^2$) hence $\pi_1(P,p)$ surjects to $\pi_1(S,w)$, with $N$ in the kernel.

A quick way to think of this is to observe that if one removes, not the entire tubular neighborhood of a surface $G$ in a 4-manifold $P$, but instead all but a single meridian disk to $G$, then the fundamental group is unchanged.  The symplectic sum of $P$ with $W_1$ can be constructed this way, where one identifies one meridian disk in $P$ with the exceptional sphere in $W_1$.
 \bigskip

 More generally, one can replace $W_1$ and $F_1$ in this remark by any appropriate pair $W,F$. We state this formally:

 \begin{lem} \label{ontolem} Suppose the 4-manifold $W$ contains a genus 2 surface $F$ with trivialized normal bundle, and the 4-manifold $P$ contains a genus 2 surface $G$ with trivialized normal bundle. Let $\phi:F\to G$ be a diffeomorphism, and let  $\tilde\phi=\phi\times Id:F\times S^1\to G\times S^1$.

 Suppose that
 \begin{enumerate}
\item  $F$ meets a sphere in $W$ transversally in one point,
\item The inclusion $F\to W$ induces a surjection on fundamental groups.
\end{enumerate}
Let $$S= (W-nbd(F))\cup_{\tilde\phi}(P-nbd(G)).$$

Then there is a surjection
$$\pi_1(P)\to \pi_1(S)$$
   whose kernel contains any word of the form $\phi(r)$, where $r\in \ker \pi_1(F)\to\pi_1(W)$.\qed
   \end{lem}

  The description preceding Lemma \ref{ontolem} indicates how to choose representative loops and base points, but in our applications of this  lemma we will typically use it  to show $S$ is simply connected, or use it when $P$ is simply connected. In either of these cases   base point issues will not matter. Notice also that choice of trivializations of the normal bundle do not affect the conclusion.

  \bigskip

  We state a similar  but easier fact whose proof can be safely left to the reader.
   \begin{lem} \label{ontolem2} Suppose the 4-manifold $W$ contains a genus 2 surface $F$ with trivialized normal bundle, with $W-F$ simply connected. Let $P$ be a 4-manifold containing a genus 2 surface $G$ with trivialized normal bundle. Let $\phi:F\to G$ be a diffeomorphism, and let  $\tilde\phi=\phi\times Id:F\times S^1\to G\times S^1$. Let
   $$S= (W-nbd(F))\cup_{\tilde\phi}(P-nbd(G)) $$

Then there is a surjection
$$\pi_1(P)\to \pi_1(S)$$
  whose kernel contains  the image of $ \pi_1(F)\to\pi_1(P)$.\qed\end{lem}

  \bigskip

\subsection{}
Our second building block $W_2$ is similar to $W_1$  but starts with $T^4$ instead  of $T^2\times S^2$:
$$W_2=T^2\times T^2\#2\bCP^2.$$  We use the calculations of Section \ref{HxK}
 to identify two Lagrangian tori $T_1$ and $T_2$ in $W_2$ and calculate the fundamental group of $W_2-(T_1\cup T_2)$, as well as their meridians and Lagrangian push offs.

Recall  from Section \ref{HxK}  that $\hat{H},\hat K$ are 2-tori,  $H$ is the complement of a small disk in $\hat H$, and $K$ is the complement of a small disc in $\hat{K}=T^2$.

Denote  $T^4\# 2\bCP^2=(\hat{H}\times \hat{K})\#2 \bCP^2$  by  $W_2$.  Then $W_2$ contains a symplectic surface $F_2$ of genus 2 with  trivial normal bundle. The construction is similar to that of $F_1\subset W_1$. Start with the symplectic surface $\hat{H}\times \{k\}\cup \{h\}\times \hat{K}\subset T^4$. Symplectically resolve the double point to obtain a symplectic surface $F_2\subset T^4$ of square $([H]+[K])^2=2$. Blow up at two points on $F_2$ to obtain $W_2$ and denote again by $F_2\subset W_2$ the proper transform.

Notice that
$W_2-nbd(F_2)$ contains the two Lagrangian tori $T_1=X\times A_1$ and $T_2=Y\times A_2$ from Theorem~\ref{tough}.   These Lagrangian tori miss the two exceptional spheres, since $\hat H \times \hat K$ is blown up at points on $F_2$, which misses $T_1$ and $T_2$.

 Recall that
$\pi_1(H\times K-(T_1\cup T_2), (h,k))$ is generated by four loops,  denoted by $x,y,a,b$ in Section \ref{HxK}.   The loops $x,y$  lie on $H\times\{k\}$ and form a basis of $\pi_1(\hat{H})$ and the loops $a,b$ lie on $\{h\}\times K$ and form a basis of $\pi_1(\hat{K})$.
Choose a small 4-ball neighborhood $B_{(h,k)}$ of $(h,k)$. Since $F_2$ is constructed by desingularizing $\hat{H}\times\{k\}\cup  \{h\}\times \hat{K}$,  we may assume that $F_2$ coincides with $\hat H\times\{k\}\cup \{h\}\times \hat K$ outside $B_{(h,k)}$.   One can choose loops $s_1,t_1,s_2,t_2$ on $F_2$  based at  point $w$ in $F_2\cap B_{(h,k)}$ which form the standard generators of $\pi_1(F_2,w)$ (in particular the relation $[s_1,t_1][s_2,t_2]=1$ holds) and which coincide with the loops  $x,y,a,b$ outside a small ball neighborhood of $(h,k)$.

\begin{prop} \label{W2}  The symplectic surface $F_2\subset W_2$ intersects an embedded sphere   transversally in one point. This sphere is disjoint from $T_1\cup T_2$, and hence the inclusion $W_2-nbd(F_2\cup T_1\cup T_2)\subset W_2-nbd( T_1\cup T_2)$ induces an isomorphism on fundamental groups, as does $Q-nbd(F)\subset Q$ for any manifold $Q$ obtained by any torus surgeries on $T_1$ and $T_2$ in $W_2$.

The fundamental group $\pi_1(W_2-nbd(T_1\cup T_2),w)$ is generated by the  loops $s_1,t_1, s_2,t_2$, which lie on $F_2$.
The relations

$$[s_1,s_2]=1, [t_1,s_2]=1, [t_1, t_2s_2t_2^{-1}]=1$$
as well as
$$[s_1,t_1]=1, [s_2,t_2]=1$$
hold in $\pi_1(W_2-nbd( T_1\cup T_2),w)$.

 Moreover, one can choose paths in $W_2-nbd(T_1\cup T_2)$ from $w$ to the boundary of the tubular neighborhoods of $T_1$ and $T_2$ so that the meridian and two Lagrangian push offs of
 $T_1$ are  expressed in $\pi_1(W_2-nbd(T_1\cup T_2),w)$ as
 $$\mu_1=[t_2^{-1}, t_1^{-1}], m_1=s_1, \ell_1=s_2$$
 and of $T_2$ are
 $$\mu_2=[s_1^{-1}, t_2], m_2=t_1, \ell_2=t_2s_2t_2^{-1}=s_2$$

\end{prop}

\begin{proof} The assertions all follow from the construction and Theorem \ref{tough}, except the relations $[s_1,t_1]=1$ and $[s_2,t_2]=1$. The relation $[s_1,t_1]=1$ holds in $W_2-(T_1\cup T_2)$ because the loops $s_1$ and $t_1$  agree with the two generators $x,y$ of the fundamental group of the torus $\hat H\times \{k\}\subset W_2-(T_1\cup T_2)$ outside a small 4-ball neighborhood of the base point.  The relation $[s_2,t_2]=1$ follows from the surface relation $[s_1,t_1][s_2,t_2]$, since $F_2$ lies in  $W_2-(T_1\cup T_2)$.

\end{proof}

Since $F_2$ meets a sphere in $W_2-nbd(T_1\cup T_2)$ transversally in one point, and $\pi_1(F_2)\to \pi_1(W_2-nbd(T_1\cup T_2))$ is surjective, Lemma \ref{ontolem} applies to the pair $(W_2-nbd(T_1\cup T_2),F_2)$. Thus if $P$ is any manifold containing a genus 2 surface $G$ with trivialized normal bundle, and $\phi:F_2\to G$ is a diffeomorphism, then the sum
$$S=(W_2-nbd(F_2))\cup_{\tilde \phi}(P-nbd(G))$$
has fundamental group a quotient of $\pi_1(P,\phi(w))$, as does
$$S-nbd(T_1\cup T_2)=(W_2-nbd(F_2\cup T_1\cup T_2))\cup_{\tilde \phi}(P-nbd(G))$$

Applying Proposition \ref{W2} we conclude that
 \begin{enumerate}
\item The kernel of the surjection $\pi_1(P,\phi(w))\to \pi_1(S-nbd(T_1\cup T_2),w)$ contains
the classes
$$\phi([s_1,s_2]), \ \phi( [t_1,s_2]),\ \phi( [t_1, t_2s_2t_2^{-1}]),\  \phi([s_1,t_1]), \ \phi( [s_2,t_2]).$$

\item The meridians and Lagrangian push offs of  $T_1$ and $T_2$ in $S$ with respect to appropriate paths to the boundary of their tubular neighborhood are given by  the images of
$$\mu_1=\phi([t_2^{-1}, t_1^{-1}]), m_1=\phi(s_1), \ell_1=\phi(s_2)$$
and
$$\mu_2=\phi([s_1^{-1}, t_2]), m_2=\phi(t_1), \ell_2=\phi(s_2)$$
under the surjection
 $$\pi_1(P,\phi(w))\to \pi_1(S-nbd(T_1\cup T_2), w)$$
 and hence if $S'$ is obtained from $1/k_i$ Luttinger surgery on $T_i$ along $\gamma_i=m_i^{p_i}\ell_i^{q_i}$ for $i=1,2$, then the kernel of the corresponding surjection
 $$\pi_1(P,\phi(w))\to \pi_1(S'))$$
 contains also the classes $\phi([t_2^{-1}, t_1^{-1}](s_1^{p_1}s_2^{q_1})^{k_1})$ and
 $\phi( [s_1^{-1}, t_2](t_1^{p_2}s_2 ^{q_2})^{k_2})$.
 \end{enumerate}

\bigskip

\subsection{} The final  and most complicated building block $M$ is a product $\hat{H}\times\Sigma $ of a torus $\hat{H}$ with a genus 2 surface $\Sigma$.
$$M=\hat H\times \Sigma.$$
Give $M$ the product symplectic form.  We will identify four Lagrangian tori $T_1,T_2,T_3,$ and $T_4$ and a genus two symplectic surface $F$ in $M$ which are pairwise disjoint and compute the fundamental group of $M-nbd(F\cup_{i=1}^4 T_i)$ and all meridians and Lagrangian push offs.

In contrast to $W_1$ and $W_2$, $M$ contains no  exceptional spheres since $\pi_2(M)=0$. In particular the inclusion
$$M-nbd(F\cup_{i=1}^4 T_i)\subset M-nbd( \cup_{i=1}^4 T_i)$$ does  not induce an isomorphism on fundamental groups. Thus we will have to be extremely careful when choosing generating loops and computing the fundamental groups of symplectic sums with $M$.

\medskip
Our approach is to view $M$ as the union of two copies of $\hat H\times K$ from Section \ref{HxK}.
 The main technical difficulty which arises   is that of identifying the generators of the fundamental group of the boundary of a tubular neighborhood of $F$ to the generators constructed from Theorem~\ref{tough}.  This is critical  in order to properly set up the use of the Seifert-Van Kampen theorem.
\bigskip

  Let $D$ be a disk with center $u$ in $\hat H$ and  identify the complement of $D$ with the surface $H$ of Section \ref{HxK}. Thus we have curves $x,y,X,Y$ and the point $h$ in $\hat H$ for Figure 1.
  To each point $q\in \hat H$, write
  $$\Sigma_q=\{q\}\times \Sigma.$$
  The surface $\Sigma_u$ corresponding to the  center $u$ of the disk $D$ will play a special role in the  following, so that we denote it by $F$:
  $$F=\{u\}\times\Sigma.$$
   The surfaces $\Sigma_q$ are symplectic for all $q$. Moreover, if $q$ misses $X\cup Y$ then $\Sigma_q$ misses all the $T_i$.
   Fix $h'$ in the boundary of the disk $D$ and choose an arc $\alpha$ in $\hat H$ joining $h'$ to $h$, as in Figure 2.
\medskip

Next view the genus 2 surface $\Sigma$ as the union of two copies of $K$ along their boundary, $\Sigma=K_1\cup_{\partial K_1=\partial K_2}K_2$.   Thus we have curves $a_1,b_1$ on $K_1$ and $a_2,b_2$ on $K_2$.
Choose arcs $\beta_1$ (resp. $\beta_2$) from  a point $k'$ on the circle separating $K_1$ and $K_2$ in $\Sigma$ to the intersection point $k_1$ of $a_1$ and $b_1$  (resp. $k_2$ of $a_2$ and $b_2$). Use the $\beta_i$ to define the corresponding based homotopy classes which satisfy $[a_1,b_2][a_2,b_2]=1$ in $\pi_1(K,k')$.  Choose two loops $A_1,A_2$ parallel to $a_1$ in $K_1$ and $A_3,A_4$ parallel to $a_2$ in $K_2$.

The notation is illustrated in Figure 2.

\begin{figure}[h]
\begin{center}
\small
\psfrag{y}{$y$} \psfrag{x}{$x$}
\psfrag{X}{$X$} \psfrag{h}{$h$}\psfrag{h1}{$h'$}\psfrag{K}{$K$}\psfrag{H}{$H$}
\psfrag{Y}{$Y$} \psfrag{A1}{$A_1$} \psfrag{A2}{$A_2$} \psfrag{A3}{$A_3$}\psfrag{A4}{$A_4$}\psfrag{a1}{$a_1$}\psfrag{a2}{$a_2$}\psfrag{b1}{$b_1$}
\psfrag{b2}{$b_2$}\psfrag{al}{$\alpha$}
\psfrag{be1}{$\beta_1$}\psfrag{be2}{$\beta_2$} \psfrag{k}{$k'$}\psfrag{u}{$u$}\psfrag{k1}{$k_1$}\psfrag{k2}{$k_2$}
 \includegraphics[scale=.90]{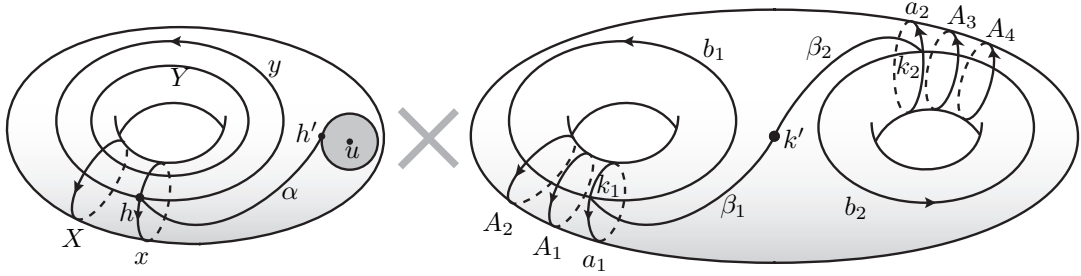}
 \caption{The surface $\hat H \times \Sigma$.}
\end{center}
\end{figure}

The product $M=\hat H\times \Sigma$ contains four disjoint Lagrangian tori $T_1=X\times A_1$, $T_2=Y\times A_2$, $T_3=X\times A_3$ and $T_4= Y\times A_4$ and  the symplectic surface
 $F=  \Sigma_u .$  These five surfaces are pairwise disjoint.

The boundary $\partial D\times \Sigma$ of  the tubular neighborhood of $F$ in $M$ contains the push off  $\Sigma_{h'}$ of $F$, as well as a meridian $\mu_F=\partial D\times \{k'\}$.  We think of
$\partial(nbd(F))$ as $\mu_F\times\Sigma_{h'} $, with base point $(h',k')$.

 The work we do in the rest of this subsection amounts to finding   loops  on $\Sigma_{h'}$ and paths between the different base points to allow us to understand the homomorphism
 $$\pi_1(  \partial (nbd(F)))\to \pi_1(M-nbd(F\cup_{i=1}^4 T_i)).$$
 explicitly.

 \bigskip

For convenience denote by $N$ the open tubular neighborhood in $M$ of the union of $F$ and the Lagrangian tori:
$$N=(D\times \Sigma)\cup nbd(T_1\cup T_2\cup T_3\cup T_4).$$
Give $M$ the base point $p=(h',k')$ on the boundary of the tubular neighborhood of $F$. We define six loops in $M-N$ based at $p$.

\begin{enumerate}
\item The loop $x\times \{k'\}$   lies on  $  H\times \{k'\}\subset M-N$  and is    based at $(h,k')$.  We conjugate this by the path $\alpha\times \{k'\}$ to define  a loop  $\tilde{x}$  based at $p=(h',k')$:
$$\tilde{x}=
 (\alpha * x*\alpha^{-1})\times \{k'\} :I\to    H\times \{k'\}\subset M-N.$$
\item The loop $y\times \{k'\}$   lies on  $ H\times \{k'\}\subset M-N$   and is    based at $(h,k')$.  We conjugate this by the path $\alpha\times \{k'\}$ to define  a loop  $\tilde{y}$  based at $p$:
$$\tilde{y}=  (\alpha * y*\alpha^{-1})\times \{k'\}:I\to    H\times \{k'\}\subset M-N.$$
\item The loops $a_1,b_1,a_2,b_2$ on $\Sigma$  defined above determine loops on
$  \Sigma_{h'}\subset \partial (D\times \Sigma)$ based at $p$:
$$\tilde a_1= \{h'\}\times a_1$$
$$\tilde b_1= \{h'\}\times b_1$$
$$\tilde a_2= \{h'\}\times a_2$$
$$\tilde b_2= \{h'\}\times b_2.$$
\end{enumerate}

Thus the loops $\tilde a_1, \tilde b_1, \tilde a_2,$ and
 $\tilde b_2$ lie on the push off $\Sigma_{h'}$ of $F$ in the  boundary of the tubular neighborhood of $F$. Together with the loop $\mu_F=\partial D\times \{k'\}$, they generate the fundamental group of $\partial(nbd(F))=\partial D \times \Sigma_{h'}$ based at $p$.

 By contrast, away from the base point,  the   loops $\tilde{x}$ and $\tilde{y}$ lie in the interior of $M-N$. However, their commutator $[\tilde x,\tilde y]$ equals $\mu_F$ in $\pi_1(M-N,p)$, since the punctured torus $H\times \{k'\}\subset M-N$  has boundary $\mu_F$.

  \bigskip

At first glance, the following proposition  may appear to be a direct application of the Seifert-Van Kampen applied to two copies of the manifold of   Theorem~\ref{tough}. However, the base point in Theorem \ref{tough} does not lie on the boundary.     Thus we must change base point {\em and} homotope appropriate loops into the boundary of $M-N$, being careful not to homotope the loops through $N$ in the process.

\begin{prop} \label{SxT}  The fundamental group  $\pi_1(M-nbd(F\cup_{i=1}^4 T_i), p)$ is generated by $\tx,\ty, \ta_1,\tb_1,\ta_2,\tb_2$ and the relations
$$1=[\tx,\ta_1]= [\ty,\ta_1]=[\ty, \tb_1\ta_1\tb_1^{-1}]=[\tx,\ta_2]=[\ty,\ta_2]=[\ty, \tb_2\ta_2\tb_2^{-1}]$$
hold in this group. With respect to certain paths to the boundary of the tubular neighborhoods of the $T_i$, the meridian and two Lagrangian push offs   are given by
\begin{enumerate}
\item $T_1: \mu_1= [\tb_1^{-1},\ty^{-1}], m_1=\tx, \ell_1=\ta_1$,
\item $T_2: \mu_2= [\tx^{-1}, \tb_1], m_2=\ty, \ell_2=\tb_1\ta_1\tb_1^{-1}$,
\item $T_3: \mu_3= [\tb_2^{-1},\ty^{-1}], m_3=\tx, \ell_3=\ta_2$,
\item $T_4: \mu_4= [\tx^{-1}, \tb_2], m_4=\ty, \ell_4=\tb_2\ta_2\tb_2^{-1}$.
\end{enumerate}
The loops $\ta_1,\tb_1,\ta_2,\tb_2$ lie on the genus 2 surface $ \Sigma_{h'} $  and form a standard set of generators  (so $[\ta_1,\tb_1][\ta_2,\tb_2]=1$).   These four loops and a meridian $\mu_F$ generate the fundamental group of the boundary of the tubular neighborhood of $F$, and $\mu_F$ is homotopic to $[\tx,\ty]$ in $\pi_1(M-nbd(F\cup_{i=1}^4 T_i),p)$.
\end{prop}

 (Please see the remark which follows the proof.)

\begin{proof}
First notice that  the punctured torus $H\times \{k'\}$ misses the tubular neighborhood $N$. Since the path $\alpha\times\{k'\}$ lies in $H\times \{k'\}$,  the boundary of this punctured torus represents the same based homotopy class as $[\tx,\ty]$ in $\pi_1(M-N,p)$. This represents the meridian $\mu_F$.

  The boundary of the tubular neighborhood of $F$ is trivialized by the push off $\Sigma_{h'}$. The curves $\ta_1,\tb_1,\ta_2,\tb_2$ lie on this  push off and so these four loops and $\mu_F$ generate the fundamental group of the boundary of the tubular neighborhood of $F$, based at $p$.

Let $S\subset \Sigma$ denote the circle separating $\Sigma$ into the two punctured tori $K_1$ and $K_2$.
Cutting $M-nbd(F)$ along $H\times S$  exhibits $M-nbd(F)$ as the union of two copies of  $H\times K$, where $H$ and $K$ are punctured tori.   The first copy $H_1\times K_1$ contains the two Lagrangian tori $T_1$ and $T_2$ and the other  contains the tori $T_3$ and $T_4$.

After cutting $M-nbd(F)$, the  surface $H\times \{k'\}$ appears as the codimension 0 submanifold
$H_1\times \{k_1'\}$ of $\partial (H_1\times K_1)$ and as the submanifold $H_2\times \{k_2'\}$ of
$\partial  (H_2\times K_2)$.  Call the copies of $\tx$ and $\ty$ that appear in $H_1\times \{k_1'\}$ $\tx_1$ and $\ty_1$, and in the other component $\tx_2$ and $\ty_2$.  The copy of $\{h'\}\times S$ (oriented and based) in $H_1\times K_1$ represents $[\ta_1,\tb_1]$ and in $H_2\times K_2$  represents $[\ta_2,\tb_2]^{-1}.$

The Seifert-Van Kampen theorem shows that
 $\pi_1(M-nbd(F\cup_{i=1}^4T_i),p)$ is the quotient of the free product
 $$\pi_1(H_1\times K_1-(T_1\cup T_2),(h_1',k_1'))*\pi_1(H_2\times K_2-(T_2\cup T_4),(h_2',k_2'))$$
 where we take the quotient by the normal subgroup  generated by $\tx_1\tx_2^{-1}$,  $\ty_1\ty_2^{-1}$
 and $[\ta_1,\tb_1][\ta_2,\tb_2]$. In particular, the loops $\tx,\ty, \ta_1,\tb_1,\ta_2,\tb_2$ generate
 $\pi_1(M-nbd(F\cup_{i=1}^4T_i),p)$.

\medskip
We reduce the proof to Theorem \ref{tough} by working one side at a time, and so, to ease   eye strain, we   drop the subscripts $1,2$. Here is what is to be shown:
We have loops  $\tx,\ty,\ta,\tb$ in $H\times K- (T_1\cup T_2)$ based at $(h',k')$ defined earlier in this section, and loops $x,y,a,b$ based at $(h,k)$ defined in the paragraph preceding the statement of Theorem \ref{tough}. The loops $x,y,a,b$ satisfy the conclusions which we will show the $\tx,\ty,\ta,\tb$   satisfy.

\medskip

 We first move from $p=(h',k')$ to $(h,k')$.   Recall we have the path $\alpha$ from $h'$ to $h$ in $H$. We let $\tilde{\alpha}$ denote the path $\alpha\times \{k'\}$.  Then conjugation by the path
 $\tilde{\alpha}^{-1}$
 defines an isomorphism
 $$\Psi_1:\pi_1(H\times K-(T_1\cup T_2),p)\to \pi_1(H\times K-(T_1\cup T_2),(h,k')), \
 \gamma\mapsto \tilde{\alpha}^{-1}*\gamma*\tilde{\alpha}.$$

From the definition preceding the statement of Proposition \ref {SxT} we see  that $\Psi_1(\tilde x)$ and $\Psi_1(\tilde y)$ are homotopic rel $(h,k')$ to the loops $x\times \{k'\}$ and $y\times \{k'\}$, since e.g.
 $$\Psi_1(\tilde x) = \tilde{\alpha}^{-1}*\tilde x*\tilde{\alpha}=
\tilde{\alpha} ^{-1}\tilde{\alpha}*(x\times\{k'\})*\tilde{\alpha}^{-1} \tilde{\alpha}\sim  x\times\{k'\}.$$

 Recall that $\ta$ takes the form $\{h'\}\times (\beta*a*\beta^{-1})$, where $\beta$ is the given  path in $K$ from $k'$ to $k$, and similarly for $\tb$.

The free homotopy $t\mapsto \{\alpha(t)\} \times (\beta*a*\beta^{-1})$  from $\ta$ to
$\{h\}\times (\beta*a*\beta^{-1}) $ misses $T_1\cup T_2$ and drags the base point along $\tilde{\alpha}$.  Hence  $\Psi_1(\tilde a)$  is represented by the loop $\{h\}\times (\beta*a*\beta^{-1})$ which lies on $\{h\}\times K$.  Similarly  $\Psi_1(\tilde b)$  is represented by the loop $\{h\}\times (\beta*b*\beta^{-1})$.

\medskip
Now we use conjugation by the path $\tilde{\beta}= \{h\}\times \beta $  to define an isomorphism
$$\Psi_2:\pi_1(H\times K-(T_1\cup T_2),(h,k'))\to \pi_1(H\times K-(T_1\cup T_2),(h,k)), \ \gamma\mapsto \tilde \beta^{-1}*\gamma*\tilde \beta$$
This takes the loop $\Psi_1(\tilde a) = \{h\}\times (\beta*a*\beta^{-1})$ to $\{h\}\times a$:
$$\Psi_2(\Psi_1(\tilde a)) =\tilde \beta^{-1}*(\{h\}\times (\beta*a*\beta^{-1}))*\tilde \beta\sim\{h\}\times  a.$$
Similarly  $\Psi_2(\Psi_1(\tilde b))=\{h\}\times b$. These are the loops simply denoted by $a$ and $b$ in Theorem \ref{tough}.

The free homotopy $t\mapsto x\times \{\beta(t)\}$ starts at $x\times \{k'\}=\Psi_1(\tilde x)$ and ends at
$x\times \{k\}$, which is the loop labeled by $x$ in Theorem \ref{tough}.  Moreover, the loop $x\times
\{\beta(t)\}$ misses $T_1\cup T_2$, since $\beta$ avoids $A_1$ and $A_2$. Since this free homotopy drags the base point along $\{h\}\times \beta=\tilde \beta$, it shows that $\Psi_2(\Psi_1(\tilde x))\sim x \times \{k\}.$ Similarly $\Psi_2(\Psi_1(\tilde y))\sim y \times \{k\}.$

\medskip

Thus we have found a path $\tau=\tilde \beta *\tilde \alpha $ in $H\times K-(T_1\cup T_2)$ from $(h,k)$ to $(h',k')$ and proven that the isomorphism
$$\pi_1(H\times K-(T_1\cup T_2),(h',k'))\to \pi_1(H\times K-(T_1\cup T_2),(h,k))$$ given by conjugating by $\tau^{-1}$  takes (the based homotopy classes of) $\tx,\ty,\ta,\tb$ to (the based homotopy classes of) $x,y,a,b$. Hence any relation satisfied by $x,y,a,b$ in $\pi_1(H\times K-(T_1\cup T_2),(h,k))$ is also satisfied by $\tx,\ty,\ta,\tb$ in $\pi_1(H\times K-(T_1\cup T_2),(h',k'))$.

  Moreover,  if one takes the paths from $(h',k')$ to the boundary of the tubular neighborhood of $T_i$ to be the composite of $\tau$ and the path given in Theorem~\ref{tough}, then e.g. the meridian of $T_1$ with respect to this path is
  $$\tau *\mu_1*\tau^{-1}= (\Psi_2\circ \Psi_1)^{-1}(\mu_1)=(\Psi_2\circ \Psi_1)^{-1}([b^{-1},y^{-1}])
  =[\tb^{-1},\ty^{-1}].$$
 A similar argument establishes the  calculations for  the other meridian and the Lagrangian push offs.

Applying the argument on each half $H_i\times K_i$ $i=1,2$ and using the Seifert-Van Kampen theorem finishes the proof.
 \end{proof}

\bigskip
\noindent{\bf Remark.} To simplify notation, for the rest of this paper  we drop the decorations, and so  we will denote $\tx$ simply by $x$ and similarly for the others. Thus the explicit loops in $M-N$ based at $p=(h',k')$ defined prior to Proposition \ref{SxT} will be denoted by $x,y,a_1,b_1,a_2,b_2$.
$a_1, b_1, a_2, b_2$ are loops that lie on $\Sigma_{h'}$ and together with  $\mu_F $ generate the fundamental group of the boundary of the tubular neighborhood of $F$.

The loops $x,y$ lie on the surface $H\times \{k'\}$ (and in particular in the interior of $M-N$ away from $p$).  The meridian $\mu_F$ equals  $[x,y]$ in $\pi_1(M-N,p)$, and the loops $x,y,a_1,b_1,a_2,b_2$ generate $\pi_1(M-N,p)$, with relations, meridians, and Lagrangian push offs as given in Proposition \ref{SxT}.

\section{Constructions of small simply symplectic manifolds}

\subsection{} We start, as a warm up, with a construction of a minimal symplectic manifold homeomorphic but not diffeomorphic to $\CP^2\#7\bCP^2$. Such examples are known
\cite{park1, OzS}; we include it  because our construction   illustrates the kind of fundamental group calculations we will do below in a simple case.

\begin{thm}\label{seven} One can perform two Luttinger surgeries on the symplectic sum of $W_1$ and $W_2$ along
$F_1$ and $F_2$ to produce a minimal symplectic manifold $U$ homeomorphic but not diffeomorphic to $\CP^2\# 7\bCP^2$. \end{thm}

\begin{proof}   Form the symplectic sum
$S=  W_1-nbd(F_1)\cup_{\tilde \phi}W_2-nbd(F_2)$ using the gluing diffeomorphism  $\phi:F_1\to F_2$
which take the loops denoted by $s_1,t_1,s_2,t_2$ on $F_1$ to their namesakes on $F_2$.
The Lagrangian tori $T_1, T_2$ in $W_2$ remain Lagrangian in $S$ (\cite[Theorem 10.2.1]{GS}).

Lemma \ref{ontolem} shows that $\pi_1(S-(T_1\cup T_2))$ is a quotient of $\pi_1(W_2-(T_1\cup T_2))$ and the kernel of the surjection contains the classes $s_1s_2, t_1t_2, $ and $[s_1,t_1]$. Applying Proposition \ref{W2} we see that
$\pi_1(S-(T_1\cup T_2))$ is a quotient of the group generated by $s_1,t_1$ and the relation $[s_1,t_1]=1$ holds, i.e. $\pi_1(S-(T_1\cup T_2))$ is a quotient of $\ZZ s_1\oplus \ZZ t_1$. Moreover, the meridians and Lagrangian push offs of the tori $T_i$ are given by
$$\mu_1=1, m_1= s_1, \ell_1=s_1^{-1}$$
and
$$\mu_2=1, m_2=t_1, \ell_2=s_1^{-1}.$$

We perform Luttinger surgeries on $T_1$ and $T_2$ in $S$  or, equivalently, Luttinger fillings on $S-nbd(T_1\cup T_2)$. Then $-1/1$ Luttinger surgery on $T_1$ along $m_1$ kills $s_1$ and $-1/1$ Luttinger surgery on $T_2$ along $m_2$ then kills $t_1$, yielding a simply connected symplectic manifold $U$.
We have $e(U)=e(W_1)+ e(W_2)+4= 10$ and $\sigma(U)=\sigma(W_1)+\sigma(W_2)= -6$. Freedman's theorem \cite{Freedman} then implies that $U$ is homeomorphic to $\CP^2\# 7\bCP^2$.

We showed that every  $-1$ sphere in $W_1$ meets $F_1$ in Proposition \ref{W1}.   Let $W_2'$ be the manifold obtained from  $W_2$ by performing the Luttinger surgeries as described, so that $U$ is the symplectic sum of $W_1$ and $W_2'$.  To see that every $-1$ sphere in $W_2'$ intersects $F_2$ takes a bit more work. Notice that $W_2'$ is obtained by performing the two Luttinger surgeries on $T_1$ and $T_2$ in $T^4=\hat H\times\hat K$  and then blowing up twice along $F_2$. These are Luttinger surgeries along $m_1$ and $m_2$, hence by Lemma \ref{aspherical},  $W_2'$ is obtained by blowing up an aspherical manifold $L$ twice along points on $F_2$.

Since $L$ is aspherical the Hopf sequence
$$\pi_2(W_2')\to H_2(W_2')\to H_2(L)\to 0$$
is exact, where $W_2'\to L$ is the map that collapses the two exceptional spheres (i.e. the blow-down map).  The kernel of $H_2(W_2')\to H_2(L)$ is clearly generated by the two exceptional spheres $E_1$ and $E_2$, and therefore every spherical class in $H_2(W_2')$ has the form $aE_1+bE_2$.   In particular, the only $-1$ spheres are $\pm E_1$ and $\pm E_2$, and both of these intersect $F$.

  If $W_1$ were an $S^2$ bundle  with section $F_1$, then $\pi_1(W_1)$ could not be $\ZZ^2$.  If $W_2'$ were an $S^2$ bundle with section $F_2$,  then the  exact sequence in homotopy groups of a fibration  would show that $\pi_2(W_2')$ equals $\pi_2(S^2)=\ZZ$. But we showed in the previous paragraph that the rank of the image of the Hurewicz map  equals 2. (One can also compute directly that $H_1(W_2')$ is generated by $s_2$ and $t_2$, so that $\pi_1(W_2')\not\cong \pi_1(F_2)$.)  Thus $W_2'$ cannot be an $S^2$ bundle over $F_2$. Applying  Usher's theorem \cite{usher}, we conclude  that $U$ is a minimal symplectic 4-manifold.

  By results of Taubes, \cite{taubes1,taubes2}, a minimal symplectic 4-manifold  cannot    contain a smoothly
embedded $-1$ sphere, but $\CP^2\#7\bCP^2$ contains smoothly embedded $-1$ spheres, namely,  the exceptional spheres. Hence $U$ cannot be diffeomorphic to  $\CP^2\#7\bCP^2$, since this would contradict the minimality of $U$.

\end{proof}

\subsection{} Our next example is more involved. We take a symplectic sum of $W_1$ and $M$ and perform four Luttinger surgeries to produce a minimal symplectic manifold homeomorphic but not diffeomorphic to $\CP^2\# 5\bCP^2$ (c.f. \cite{PSS, A}).

Consider the surface $H$  and its curves $x,y,X,Y$  and base point $h$ in Figure 1.  Let $D_X$ and $D_Y$ denote the Dehn twists along $X$ and $Y$. We leave the proof of the following simple lemma to the reader.

\begin{lem}\label{dehn} The composite of Dehn twists along $X$ and $Y$ given by $D_XD_YD_X$ takes $x$ to $y^{-1}$ and $y$ to $yxy^{-1}$ in $\pi_1(H,h)$.
\qed\end{lem}
 Note Lemma \ref{dehn} holds in $\pi_1(H,h')$ if we replace $x,y$ by $\al*x*\al^{-1}$, $\al*y*\al^{-1}$ where $\al$ is the path in $H$ from the boundary point $h'$ to $h$ which misses $X$ and $Y$. Moreover,  $D_XD_YD_X$ is the identity on the boundary, so that Lemma \ref{dehn} extends to higher genus surfaces.

 \medskip

Consider the surface $F_1\subset W_1$ with base point $w$ and $\Sigma_{h'}\subset M$ with base point $p=(h',k')$.  Denote by  $\phi:F_1\to \Sigma_{h'}$ the base point preserving diffeomorphism  of Lemma \ref{dehn}  inducing the isomorphism
$$\phi:\pi_1(F_1, w)\to \pi_1(\Sigma_{h'},p), \ (s_1,t_1,s_2,t_2)\mapsto (b_1^{-1},b_1a_1b_1^{-1},a_2,b_2)$$
and extend   to $\tilde{\phi}:F_1\times S^1\to \Sigma_{h'}\times S^1$.
The symplectic sum
$$ V'=(W_1-nbd(F_1))\cup (M-nbd(\Sigma_{h'}))$$
contains the four Lagrangian tori $T_1,T_2,T_3,T_4$  and a symplectic surface $F$ so that these five surfaces are pairwise disjoint.

Using Lemma \ref{ontolem} and the facts that $s_1=s_2^{-1}$, $t_1=t_2^{-1} $, and $[s_1,t_1]$ in $\pi_1(W_1)$, we see that $b_1=a_2$ in $\pi_1(V'-nbd(F\cup_{i=1}^4 T_i),w)$. The relation  $[s_1,t_1]=1$ implies that  $1=[b_1^{-1}, b_1a_1b_1^{-1}]=[a_1,b_1^{-1}]$ and hence $[a_1,a_2]=[a_1,b_1]=1$ and $b_2=a_1^{-1}$
hold in $\pi_1(V'-nbd(F\cup_{i=1}^4 T_i),w)$.

Therefore, $\pi_1(V'-nbd(F\cup_{i=1}^4 T_i),w)$ is generated by $x,y,a_1,a_2$ and the relation
$$[a_1,a_2]=1$$
holds.
The relations
$$[x, a_1]=1, [y, a_1]=1, [x, a_2]=1, [y, a_2]=1$$ coming from $M-nbd(F\cup_{i=1}^4 T_i)$ were established in Proposition \ref{SxT}, and by the Seifert-Van Kampen theorem also hold
in $\pi_1(V'-nbd(F\cup_{i=1}^4 T_i))$.

Furthermore, it follows from Proposition \ref{SxT} that
the meridian of $F$ is $$\mu_F=[x,y].$$
and the meridians and Lagrangian push offs of the tori  $T_1$ are given by
\begin{itemize}
\item $T_1: \mu_1= [a_2^{-1}, y^{-1}]=1, m_1= x, \ell_1= a_1$.
\item $T_2: \mu_2= [ x^{-1},  a_2]=1, m_2= y, \ell_2= a_2 a_1a_2^{-1}= a_1$.
\item $T_3: \mu_3= [a_1, y^{-1}]=1, m_3= x, \ell_3= a_2$.
\item $T_4: \mu_4= [ x^{-1}, a_1^{-1}]=1, m_4= y, \ell_4=a_1^{-1} a_2a_1= a_2$.
\end{itemize}

\begin{thm}\label{five} There exists a minimal symplectic manifold $V$ homeomorphic but not diffeomorphic to $\CP^2\# 5\bCP^2$ containing a symplectic genus 2 surface $F$ with simply connected complement and trivial normal bundle.
\end{thm}
\begin{proof}  We do four Luttinger  surgeries on $T_1, T_2, T_3, $ and $T_4$ in $V'-F$.
First, $-1/1$ Luttinger surgery on $T_1$ along $m_1$ kills $x$. Then $-1/1$ surgery  on $T_2$ along $\ell_2$ kills $a_1$. Next, $-1/1$ surgery on $T_3$ along $\ell_3$ kills $a_2$. Finally, $-1/1$ surgery on $T_4$ along $m_4$ kills $y$.

Thus we have produced a symplectic 4-manifold $V$ such that  $V-F$ is simply connected. The Seifert-Van Kampen theorem shows that $V$ is simply connected. The Euler characteristic and signature  of $V$ are the same as $V'$ since they are unchanged by Luttinger surgery. Hence
$$e(V)=e(W_1)+e(M)+4=4+0+4=8$$
and
$$\sigma(V)=\sigma(W_1)+\sigma(M)= -4+0=-4.$$
Freedman's theorem \cite{Freedman} then implies that $V$ is homeomorphic to $\CP^2\# 5\bCP^2$.

Let $M'$ be the manifold obtained from $M=\hat H\times\Sigma$ by doing the four Luttinger surgeries described above.    Minimality of $V$  follows from Usher's theorem \cite{usher} just as in the proof of Theorem \ref{seven} once we show that $M'$   is minimal, and not an $S^2$ bundle with section $\Sigma_{h'}$.  Note that $M'$ is the symplectic sum of two manifolds $M'_1$ and $M'_2$ along the torus $\hat H$.

 Indeed,  view $M=\hat H\times \Sigma$ as the union of $\hat H_1\times K_1$ and $\hat H_2\times K_2$ along $\hat H\times S$ as in the proof of Proposition \ref{SxT}.  This exhibits $M$ as the symplectic sum of $\hat H_1\times \hat K_1$ and $\hat H_2\times \hat K_2$ along  the symplectic surfaces $\hat H_1\times \{v_1\}$ and $\hat H_2\times \{v_2\}$, where $v\in \hat K$ is the center of the disk whose complement is $K$.  So let $M_1'$ denote the manifold obtained by Luttinger surgery on $T_1$ and $T_2$ in the 4-torus $\hat H_1\times \hat K_1$ and  $M_2'$ denote the manifold obtained by Luttinger surgery on $T_3$ and $T_4$ in the 4-torus $\hat H_2\times \hat K_2$.

Lemma \ref{aspherical} implies that $M_1'$ and $M_2'$ are aspherical, hence minimal. In $H_1(M_1')$, the generators $x$ and $a_1$ are zero, so that $H_1(M_1')$ is generated by $y$ and $b_1$. In particular, $\hat H_1\times \{v_1\}\subset  M_1'$ cannot be a section of an $S^2$ fiber bundle structure on $M_1'$ since $\pi_1(\hat H_1\times \{v_1\})$ is generated by $x$ and $y$. Thus Usher's theorem implies that $M'$ is minimal.  Finally, $M'$ is not an $S^2$ bundle with section $\Sigma_{h'}$ for similar reasons: the loops $a_1,b_1,a_2,b_2$ are generators of $\pi_1(\Sigma_{h'})$ but $a_1$ and $a_2$ are trivial in $H_1(M')$, as one can readily check. As explained above, it follows that $V$ is minimal.

Using the same argument as in the proof of Theorem \ref{seven} one concludes from Taubes's results   that $V$ is not diffeomorphic to $\CP^2\#5\bCP^2$.
\end{proof}

Since $V$ contains an appropriate Lagrangian torus, the argument of  Corollary \ref{infinitely} below applies to the manifold $V$ as well to produce infinitely many smooth (but not symplectic)  pairwise non-diffeomorphic manifolds homeomorphic to $V$.
\bigskip

\subsection{}  We next put together the manifolds $W_2$ and $M$.  The construction is based on the example of a  minimal symplectic manifold homeomorphic to $\CP^2\# 3\bCP^2$ which we constructed in \cite{BK2}. However, the extra information obtained by keeping track of the surface $F$ will allow us to produce many more examples which we will use in constructing small non-simply connected symplectic 4-manifolds    in the next section.

\bigskip

Recall that the manifold $W_2$ contains a genus 2 square zero symplectic surface $F_1$ and two Lagrangian tori $T_1$ and $T_2$.

The manifold $M=\hat H\times\Sigma$ of Proposition \ref{SxT} contains a symplectic surface $F$ with trivial normal bundle and four Lagrangian tori $T_1,T_2,T_3,T_4$, with the fundamental group of
$M-nbd(F\cup T_1\cup T_2\cup T_2\cup T_4)$
generated by the loops $ x,y,a_1,b_1,a_2,b_2$ based at $p$ satisfying all the conclusions of Proposition \ref{SxT}.

 To avoid notational confusion, we denote the two Lagrangian tori in $W_2$ by $T_1'$ and $T_2'$. We preserve the notation $T_1,T_2,T_3,T_4$ for the four Lagrangian tori in $M$.

 \medskip

The  parallel (and symplectic) push off $\Sigma_{h'}$ of $F$   lies in the boundary of the tubular neighborhood $D\times \Sigma$ of $F$ and carries the base point and the  loops $a_1,b_1,a_2,b_2$.      We have a framing of $\Sigma_{h'}\subset M$ defined by taking a nearby   push off $\Sigma_{z}$   for a point $z$ near $h'$. Choose some identification  of the tubular neighborhood of $F_2\subset W_2$ with $F_2\times D^2$.

\bigskip
We form the symplectic sum $Z$ of $W_2$ and $M$ along  $F_2$ and $\Sigma_{h'}$:
\begin{equation}\label{Zmanifold}
Z=W_2-nbd(F_2)\cup_{\tilde\phi}M-nbd(\Sigma_{h'}).
\end{equation}
Observe that we have taken the symplectic sum along $\Sigma_{h'}$, not $F$. Thus $F$ survives as a symplectic genus $2$ surface in $Z$.

\bigskip

We choose the base point preserving diffeomorphism
$\phi:F_2\to \Sigma_{h'}$ to form this sum as follows. Using Lemma \ref{dehn} in each half of the decomposition of $\Sigma$ into two punctured tori, we conclude that there is a base point preserving diffeomorphism (a composite of six Dehn twists) $\phi_1:\Sigma_{h'}\to\Sigma_{h'}$ that induces the isomorphism
$a_1,b_1,a_2,b_2$ to $ b_1^{-1}, b_1a_1b_1^{-1}, b_2^{-1} , b_2a_2b_2^{-1}$ on $\pi_1(\Sigma_{h'},h')$.
Composing this with the diffeomorphism $\phi_2:F_2\to \Sigma_{h'}$ that takes   $s_1,t_1,s_2,t_2$ to  $a_1,b_1,a_2,b_2$ yields the desired diffeomorphism $\phi:F_2\to \Sigma_{h'}$.  Hence $\phi$ induces the isomorphism
\begin{equation}\label{phi}
\phi: \pi_1(F_2, w)\to \pi_1(\Sigma_{h'},h'), \ \ (s_1,t_1,s_2,t_2)\mapsto (   b_1^{-1}, b_1a_1b_1^{-1}, b_2^{-1} , b_2a_2b_2^{-1}).\end{equation}

The symplectic manifold $Z$ contains the surface $F$ and six Lagrangian tori $T_1',T_2',T_1,T_2,T_3,T_4$. These seven surfaces are pairwise disjoint.  For convenience
denote the union of  these seven surfaces by $R\subset Z$.
 Lemma \ref{ontolem} and the discussion following Proposition \ref{W2} shows that
 $\pi_1(Z-R)$ is a quotient of $\pi_1(M-(F\cup_{i=1}^4 T_i))$.  Thus
$\pi_1(Z-R)$ is generated by   $x,y,a_1,b_1,a_2,b_2$.

Moreover, the relations of Proposition \ref{W2} and Lemma \ref{ontolem} imply that in $\pi_1(Z-R)$,
$$[\phi(s_1),\phi(s_2)]=[\phi(t_1),\phi(s_2)]=[\phi(s_1),\phi(t_1)]=[\phi(s_2),\phi(t_2)]=1.
$$
Rewriting this in terms of the $a_i,b_i$ using Equation (\ref{phi}) one obtains
$$
1= [b_1^{-1}, b_2^{-1} ] ,  1= [b_1a_1b_1^{-1}, b_2^{-1}], 1=[b_1^{-1} , b_1a_1b_1^{-1}] , 1= [b_2^{-1} , b_2a_2b_2^{-1}].
$$
Notice that $[b_1^{-1},b_1a_1b_1^{-1}]=[a_1,b_1]$ and $[b_2^{-1},b_2a_2b_1^{-1}]=[a_2,b_2]$.  Moreover, $[r,s]=1$ implies $[r^{-1},s^{-1}]=1$ and $[r,s^{-1}]=1$. Hence this set of relations simplifies to
\begin{equation}\label{eq4.1}
1= [b_1 , b_2  ] ,  1= [ a_1 , b_2 ], 1=[b_1  ,  a_1 ] , 1= [b_2  , a_2 ].
\end{equation}

Since $\pi_1(Z-R)$ is a quotient of $\pi_1(M-nbd(F\cup_{i=1}^4 T_i))$, the relations coming from $M-nbd(F\cup_{i=1}^4 T_i)$ also hold in $\pi_1(Z-R)$ , hence
$$
 [x,a_1]= [y,a_1]=[y, b_1a_1b_1^{-1}]=[x,a_2]=[y,a_2]=[y, b_2a_2b_2^{-1}]=1.
$$
These simplify, using Equation (\ref{eq4.1}),  to
\begin{equation}\label{eq4.2}
1=[x,a_1], 1= [y,a_1], 1 =[x,a_2], 1=[y,a_2] .
\end{equation}

\bigskip

The following is the result we have been aiming towards. Its usefulness will be illustrated in most of the subsequent   constructions in this article.

 \begin{thm}\label{prop1} Let $Z$ denote the symplectic sum of $W_2$ and $M$ along the surfaces $F_2$ and $\Sigma_{h'}$ using the diffeomorphism $\phi:F_2\to \Sigma_{h'}$ inducing the isomorphism of  Equation (\ref{phi}) on fundamental groups. Let $R$ denote the union the seven surfaces in $Z$: the symplectic surface $F$ and the six Lagrangian tori $T_1',T_2',T_1,T_2,T_3,T_4$.

 Then the fundamental group of $Z-R$ is  generated by loops $x,y,a_1,b_1,a_2,b_2$. The relations
 $$1= [b_1 , b_2  ] ,  1= [ a_1 , b_2 ], 1=[b_1  ,  a_1 ] , 1= [b_2  , a_2 ] $$
and
$$1=[x,a_1], 1= [y,a_1], 1 =[x,a_2], 1=[y,a_2] $$
 hold in $\pi_1(Z-R)$.

 The meridian of the surface $F$ is given by $\mu_F=[x,y]$. For the six Lagrangian tori, the meridians and Lagrangian push offs (with appropriate paths to the base point) in $\pi_1(Z-R)$ are given by:
\begin{itemize}
\item $T_1': \mu_1' =[a_2^{-1}, a_1^{-1}], m_1'=b_1^{-1}, \ell_1'=b_2^{-1}$.
\item  $T_2': \mu_2'= [b_1,a_2], m_2'= b_1a_2b_1^{-1},
\ell_2'=  b_2^{-1}$.
\item $T_1: \mu_1= [b_1^{-1},y^{-1}], m_1=x, \ell_1=a_1$.
\item $T_2: \mu_2= [x^{-1}, b_1], m_2=y, \ell_2  =a_1$.
\item $T_3: \mu_3= [b_2^{-1},y^{-1}], m_3=x, \ell_3=a_2$.
\item $T_4: \mu_4= [x^{-1}, b_2], m_4=y, \ell_4 =a_2$.
\end{itemize}

  \end{thm}

  \begin{proof}    Proposition \ref{SxT} implies that the meridian of $F$ in
 $\pi_1(Z-R)$ is given by
$$
 \mu_F=[x,y].
$$

 Propositions \ref{W2} and \ref{SxT} give the meridians and Lagrangian push offs of the six tori in terms of these generators.

  Proposition \ref{W2} shows that
  $\mu_1'=[b_2a_2^{-1}b_2^{-1}, b_1a_1^{-1}b_1^{-1}]$.
 Equation  (\ref{eq4.1})  can be used to simplify this:  $[b_2a_2^{-1}b_2^{-1}, b_1a_1^{-1}b_1^{-1}]=[a_2^{-1}, a_1^{-1}]$, and hence $\mu_1'$ is equal to $[a_2^{-1}, a_1^{-1}]$.  Similarly $[b_1, b_2a_2b_2^{-1}]=[b_1,a_2]$, showing that $\mu_2'=[b_1,a_2]$. Next, $(b_2a_2b_2^{-1})b_2^{-1}(b_2a_2b_2^{-1})^{-1}=b_2^{-1}$, showing that $\ell_2'=  b_2^{-1}$.
Continuing, $ b_1a_1b_1^{-1}=a_1$ and so $\ell_2=a_1$ and  $b_2a_2b_2^{-1}=a_2$, so that
  $ \ell_4= a_2$.

  The rest of the argument was described before the statement of the Theorem.

  \end{proof}

  We determine the basic homological properties of $Z$.

  \begin{prop}\label{sixtori} The first homology $H_1(Z)\cong \ZZ^6$, generated by the loops $x,y,a_1,b_1, a_2,b_2$.

  The second homology
  $H_2(Z)\cong \ZZ^{16}$.
There exist six disjoint tori $R_1',R_2', R_1,R_2,R_3,$  and $R_4$, and a genus 2 surface $H_3$ with trivial normal bundles  in $Z$ which are geometrically dual to the   six Lagrangian tori  $T_1',T_2',T_1, T_2,T_3$, and $T_4$ and the surface $F$, in the sense that $T'_i$ intersects $R'_i$ transversally in one point and similarly for $T_i, R_i$ and   $F,H_3$, and all other intersections are pairwise empty.

    There are two disjoint tori $H_1$, $H_2$ with square $-1$ which intersect $F$ transversally once and are pairwise disjoint from the $T_i',R_i', T_i,R_i,$ and $H_3$. The sixteen surfaces $T_i',R_i',T_i,R_i, F,$ and $H_1$ generate $H_2(Z)$. Hence the intersection form of $Z$ is an orthogonal sum of six hyperbolic summands generated by the pairs $T_i,R_i$ (and $T_i',R_i')$ and one 4-dimensional summand spanned by $H_1,H_2,H_3$ and $F$ with matrix
$$\begin{pmatrix} -1&0&0&1\\
0&-1&0&1\\
0&0&0&1\\
1&1&1&0
\end{pmatrix}$$

  The signature $\sigma(Z)$ equals $-2$. The Euler characteristic $e(Z)$ equals $6$.

  Moreover, the symplectic manifold $Q$ obtained by performing any six Luttinger surgeries on the $T_i'$ and $T_i$  satisfies $e(Q)=6$ and $\sigma(Q)=-2$, and contains a symplectic genus 2 surface $F$ with trivial normal bundle and six Lagrangian tori.

 \end{prop}

 \begin{proof} Since $Z$ is the symplectic sum of $W_2=T^4\#2\bCP^2$ and $M=T^2\times \Sigma$ along  genus 2 surfaces,
 $$e(Z)=e(T^4\#2\bCP^2)+e(T^2\times \Sigma) +4=6$$
 and Novikov additivity implies that
 $$\sigma(Z)=\sigma(T^4\#2\bCP^2)+\sigma(T^2\times \Sigma) =-2.$$

 Luttinger surgery does not affect $e$ and $\sigma$, and so   $e(Q)=6$ and $\sigma(Q)=-2$.  The core torus $T\times \{0\}\subset T^2\times D^2$ glued in to form a Luttinger surgery is Lagrangian \cite{Lut}. Since the $T_i'$ and $T_i$ miss $F$,  the assertions about $Q$ are verified.

 The  homomorphism $H_1(F_2)\to H_1(W_2-F_2)=H_1(W_2)\cong \ZZ^4$ induced by any push off is an isomorphism. Also $M-nbd(\Sigma_{h'})= H\times \Sigma$, where $H$ is a punctured torus. Thus from the Mayer-Vietoris sequence for the decomposition
$$Z=W_2-nbd(F_2)\cup_{\tilde\phi}M-nbd(\Sigma_{h'})$$
one easily sees that $H_1( M-nbd(\Sigma_{h'}))\to H_1(Z)$ is an isomorphism.
Since $H_1(M-nbd(\Sigma_{h'}))=H_1(H\times \Sigma)\cong \ZZ^6$, and $\pi_1(H\times \Sigma)$ is generated by the six loops $x,y,a_1,b_1,a_2,b_2$, $H_1(Z)\cong\ZZ^6$ generated by these loops.
Since $e(Z)=6$, it follows that $H_2(Z)\cong \ZZ^{16}$.

\medskip

We find 16 embedded surfaces that generate $H_2(Z)$ and calculate the intersection form.

The torus $T_1$ has a geometrically dual torus $R_1=x\times a_1$ in $M=\hat H\times \Sigma$ which misses $T_2,T_3,T_4$ and also $F$ and its parallel $\Sigma_{h'}$ (since the points $u$ and $h'$ are different from $h$).   Similarly $T_2$ has a dual torus $R_2=x\times b_1$, $T_3$ has a dual torus $R_3=y\times b_2$, and $T_4$ has a dual torus $R_4=x\times b_2$.
 These  all miss $F$ and $\Sigma_{h'}$ and so survive in $Z$. It is also easy to push the $R_i$ off each other. Notice also that the $T_i$ and $R_i$ miss $\hat H\times \{q\}$ for most points $q$.

Similarly one can find disjoint dual tori $R_1'$ and $R_2'$ to the classes $T_1'$ and $T_2'$ in $W_2-F_2$.  The tori $T_i'$ and $R_i'$ may be assumed to miss the 2 exceptional curves $E_1$ and $E_2$, and also miss one of the vertical tori $\{r\}\times T^2$.

Let $H_1\subset Z$ be the torus formed by joining up the exceptional curve $E_1$ to   $H\times\{q\}$
for the appropriate  $q$ in the symplectic sum. Similarly define $H_2$ using the other exceptional curve, and let $H_3$  be a genus 2 surface formed by joining $H\times\{r\}$ to one of the vertical tori.

 The  12 tori $T_i, R_i, T_i', R_i'$ can be isotoped by a small isotopy in $Z$ so  each meets its dual transversally once and all other intersections are empty.   Moreover, they each miss $F$, $H_1$, $H_2$, and $ H_3$.  Notice that each $H_i$ intersects $F$ once (geometrically and algebraically), that $F^2=0$, $H_1^2=-1$,   $H_2^2=-1$, and $H_3^2=0$. Finally $H_i$ is disjoint from $H_j$ when $i\ne j$.   One can check that the classes $F=H_0,H_1,H_2,H_3$ span a primitive subspace of $H_2(Z)$ by calculating that the determinant of $H_i\cdot H_j$ is equal to $1$.

Thus the intersection form of $Z$ is the orthogonal sum of 6 hyperbolic planes spanned by each $T_i,R_i$ pair, and a 4 dimensional space spanned by $F,H_1,H_2,H_3$ with matrix as asserted in the statement.

 Notice moreover that  the subspace spanned by $F$ and the $T_i$ is a 7 dimensional isotropic subspace,  and that $H_2$ and $H_3$ span a subspace with intersection form $2(-1)$.

 \end{proof}

  \bigskip

We can now do Luttinger surgery on the six tori to obtain  interesting symplectic 4-manifolds containing a symplectic genus 2 surface $F$.  The effect of Luttinger surgery on the homology is easy to understand.
Let  $T$ denote  one of the $T_i$ or $T_i'$, and $\mu,m,\ell$ its meridian and two Lagrangian push offs.
Theorem \ref{prop1} shows  that $m$ is sent to one of the generators $x,y, a_1,b_1,a_2,b_2$ in $H_1(Z)$, and similarly for $\ell$. The meridian $\mu$ is trivial in $H_1(Z)$.  Therefore,

\begin{enumerate}
\item   $1/0$ surgery on $T$ along any curve $\gamma$ does not change anything, one removes a neighborhood of $T$ and re-glues it the same way.
\item For $k\ne 0$, $1/k$ surgery on $T$ along the curve $\gamma= m^p\ell^q$ kills $k(p[m]+q[\ell])$ in $H_1(Z)$. The surgery decreases the rank of $H_1(Z;\QQ)$ by one.  Since the Euler characteristic is unchanged, the rank of $H_2(Z;\QQ)$ decreases by two. If $k=\pm 1$, then the rank of $H_2(Z;\QQ)$ decreases by two.
\end{enumerate}

Notice that $\pm 1/1$ Luttinger surgery  on  $T$ only changes the manifold near $T$. Letting $\tilde{T}$ denote the core $T\times \{0\}$ of $T\times D^2$ in the {\em surgered} manifold,  we see that $\tilde{T}$ is nullhomologous, since the other surfaces $T_i$, $R_i$, $H_i$ and  $F$ are disjoint from $T$ and its dual and hence are not affected. In other words, the effect on second homology of $\pm 1$ Luttinger surgery on $T_i$ or $T_i'$ is to kill the hyperbolic summand spanned by $T_i$ and its dual.

\medskip

 Computing fundamental groups is harder, but the following refinement of the main result of  \cite{BK2} indicates why
we chose the gluing map $\phi$ as we did.

\begin{thm}\label{cool} There exists a minimal symplectic 4-manifold $X$ homeomorphic but not diffeomorphic to $\CP^2\# 3\bCP^2$ containing a symplectic surface $F$ of genus 2    with simply connected complement and trivial normal bundle.

Moreover, $X$ contains a nullhomologous Lagrangian torus $T$ disjoint from $F$, and a curve $\lambda$ on the boundary of the tubular neighborhood of $T$ such that $\lambda=0$ in $H_1(X-T)$ and $\lambda$ is isotopic in $nbd(T)$ to an embedded essential curve  on $T$.

\end{thm}
\begin{proof} We perform six Luttinger surgeries on the six Lagrangian tori in $Z-F$, referring to Theorem \ref{prop1}.  First $1/1$ surgery on $T_1'$ along $m'_1$ imposes the relation $b_1=[a_2^{-1}, a_1^{-1}]$. Then $-1/1$ surgery on $T_1$ along $m_1$ imposes the relation $x=[b_1^{-1}, y^{-1}]= [[a_1^{-1}, a_2^{-1}], y^{-1}]=1$. The last equality comes from Equation
(\ref{eq4.2}).

Next perform $-1/1$ Luttinger surgery on $T_2$ along $\ell_2$. This imposes
$a_1=[x^{-1},b_1]=1$.   Returning to our first surgery we see that $b_1=[a_2^{-1}, a_1^{-1}]=1$.
Continuing, $ 1/1$ surgery on $T_2'$ along $\ell_2'$ sets
$b_2=[b_1,a_2]=1$.

Next perform $-1/1$ surgery on $T_3$ along $\ell_3$ to set $a_2=1$. We finish with  $-1/1$ surgery on $T_4$ along $m_4$ which sets $y=1$.

This kills all the generators, so the resulting symplectic 4-manifold $X$ satisfies $\pi_1(X-F)=1$.

 It was explained in the paragraphs preceding Theorem \ref{cool} that the  six Luttinger fillings kill the first homology, and $H_2(X)=\ZZ^4$,   and that any of the six core tori in $X$ are nullhomologous.
Hence we let $T$ denote the core of the neighborhood glued in the last Luttinger surgery.

 Since $T\times D^2$ is glued along $T\times\partial D^2$ so that  $\partial D^2$ is identified with the curve $\mu_4 m_4^{-1}$, we can choose an essential embedded curve $\lambda$ on $T$ so that $\lambda\times\{s\}$ is sent to $\mu_4=[x^{-1},b_2]=1$ in $\pi_1(X-T)$, and hence
 $\mu_4=0$ in $H_1(X-T)$. (Note that since each Luttinger surgery kills one generator of first homology, $\pi_1(X-T)=H_1(X-T)=\ZZ y$.)

The minimality of $X$ follows just as in the proofs of Theorems \ref{seven} and \ref{five}, using Usher's theorem.   Note that $e(X)=e(Z)=6$ and $\sigma(X)=\sigma(Z)=-2$, so that by Freedman's theorem \cite{Freedman} $X$ is homeomorphic to $\CP^2\#3\bCP^2$.  As in the proof of Theorem \ref{seven}, Taubes's results imply that $V$ is not diffeomorphic to $\CP^2\#3\bCP^2$.
This completes the proof.
 \end{proof}

Referring to Proposition \ref{sixtori}, we see that the four classes generating $H_2(X)$ and diagonalizing the intersection form are $F+H_1$, $H_1$, $H_2-H_3$, $H_3-H_1-F$. These are represented by smoothly embedded surfaces of genus $3,1,3,$ and $5$ respectively. Notice that $H_1$ has minimal genus, since $X$ cannot contain a smoothly embedded $-1$ sphere. Similarly $H_2$ has minimal genus.

\bigskip

\noindent{\bf Remark.} The last Luttinger surgery in the proof of Theorem \ref{cool}, $-1/1$ surgery on $T_4$, kills the last loop $m_4=y$.   But one can leave the setting of symplectic manifolds, and
define a family of smooth manifolds as follows. Denote by $X_0$ the symplectic manifold constructed in the proof of Theorem \ref{cool} in the penultimate step. Thus $X_0 $ is a symplectic manifold with $\pi_1(X_0)=\ZZ$, generated by $y$,  and $b^+(X_0)=2, b^-(X_0)=4$.  The torus $T_4$ in $X_0$ is Lagrangian. The boundary of the tubular neighborhood of $T_4$ is a 3-torus whose fundamental group is generated by the loops $\mu_4,m_4$ and $\ell_4$. In $\pi_1(X_0-T_4)$, $\mu_4=1$, $m_4=y$, and $\ell_4=1$.

Let $Y$ denote the manifold obtained from $X_0-nbd(T_4)$ by gluing $T^2\times D^2$ in such a way that $\alpha=S^1\times \{1\}\times \{1\}$ is sent to  $\ell_4$, $\beta=\{1\}\times S^1\times \{1\}$ is sent to $\mu_4$, and $\mu_Y=\{(1,1)\}\times \partial D^2$ is sent to $m_4^{-1}$.   Then $Y$ is simply connected, since $m_4=y$ is killed, but $Y$ need not be symplectic.   Let $\tilde T\subset Y$ denote the resulting core   torus.

With these coordinates $\alpha, \beta, \mu_Y$, one sees that
\begin{enumerate}
\item $X_0$ is obtained by $0/1$ surgery on $\tilde T$ along $\beta$,
\item $\beta$ is nullhomologous (indeed nullhomotopic) in $Y-\tilde T=X_0-T_4$, since $\beta=\mu_4$.
\item the manifold $Y_n$ obtained from $1/n$ surgery on $\tilde T$ along $\beta$ is simply connected, and hence homeomorphic to $\CP^2\#3\bCP^2$. (This manifold can be described as $-n/1$ surgery on
$T_4$ along $m_4$ in $X_0$, and so $\mu_4^nm_4^{-1}=y^{-1}$ is killed. Thus $X=Y_1.$)
\item the torus $\tilde T$ is nullhomologous in $Y$.  The reason is the same as in the construction of $X$: the last  torus surgery kills one generator of $H_1$, hence two generators of $H_2$, the class of $T_1$ and its dual $R_1$.
\item $X_0$ is symplectic, minimal, and has $b^+=2$, hence has non-trivial Seiberg-Witten invariant by \cite{taubes3}. \end{enumerate}

Thus, as explained by Fintushel and Stern in \cite{FS5} (see also \cite{FS2}), the Morgan-Mrowka-Szabo formula \cite{MMS} can be used to prove that the family $\{Y_n\}$ contains infinitely many diffeomorphism types, detected by Seiberg-Witten invariants. This proves the following.

\begin{cor} \label{infinitely} Among the manifolds $Y_n, n\in\ZZ$ there exists an infinite family of smooth manifolds homeomorphic to $\CP^2\# 3\bCP^2$ which are pairwise non-diffeomorphic (detected by Seiberg-Witten invariants). Each $Y_n$ contains a smoothly embedded genus 2 surface $F$ with trivial normal bundle and simply connected complement. \qed
\end{cor}

 As remarked in \cite{SS}, Seiberg-Witten invariants cannot be used to distinguish more than two minimal symplectic manifolds homeomorphic to $\CP^2\#n\bCP^2$  for $n<9$. Thus all but at most two of the $Y_n$ are minimal symplectic manifolds.

\bigskip

\subsection{} The manifold   $X$ is very useful. We illustrate this with a number of constructions next.
 A   quick corollary is the following (cf. \cite{A}).

\begin{cor}\label{b+=31} The symplectic sum of $X$ with $W_2$ along  $F$ in $X$ and $F_2$ in $W_2$  is a minimal symplectic manifold homeomorphic but not diffeomorphic to $3\CP^2\# 7\bCP^2$.
\end{cor}

\begin{proof} The map $\pi_1(F_2)\to \pi_2(W_2)$ induced by inclusion is onto, and hence the symplectic sum $S$ of $W_2$ and $X$ is simply connected by  Lemma \ref{ontolem2}.  Since
 $e(S)=e(X)+e(W_2)+4=12$ and $\sigma(S)=\sigma(X)+\sigma(W_2)=-4$,   Freedman's theorem \cite{Freedman} implies that $S$ is homeomorphic to $3\CP^2\# 7\bCP^2$.

The Seiberg-Witten invariants of a connected sum of manifolds with $b^+>1$ vanishes (\cite{Sal}).  On the other hand the Seiberg-Witten invariants for a symplectic manifold with $b^+>1$ are non-trivial (\cite{taubes3}).   Hence $S$ cannot be diffeomorphic to $3\CP^2\# 7\bCP^2$. \end{proof}

Similarly one obtains the following (cf. \cite{SS1}).
\begin{cor}\label{b+=32}  The fiber sum of $X$ with $W_1$ along $F$ and $F_1$ is a minimal symplectic manifold homeomorphic but not diffeomorphic to $3\CP^2\# 9\bCP^2$. \qed
\end{cor}

and also

\begin{cor}\label{b+=51}  The fiber sum of two copies of $X$   along $F$  is a minimal symplectic manifold homeomorphic but not diffeomorphic to $5\CP^2\# 9\bCP^2$. \qed
\end{cor}

In Corollaries \ref{b+=31}, \ref{b+=32}, and \ref{b+=51} the manifolds constructed contain nullhomologous Lagrangian tori  with an appropriate curve to surger, and so the argument given above for $X$ applies to find infinitely many smooth homeomorphic  but not diffeomorphic (and not necessarily symplectic) examples in each case.

\bigskip

\subsection{} It can be shown that any symplectic 4-manifold containing a symplectic torus of square zero either has $e\ge 12$ or $b^+\ge 3$.  Thus the small simply connected symplectic manifold $X$ with $b^+=1$ constructed above  cannot contain a square zero torus with simply connected complement.

 The elliptic surface $E(1)$, which has $e=12$ and $\sigma=-8$, does contain a symplectic square zero torus with simply connected complement, namely the fiber in any elliptic fibration with cusp fibers.
 The following theorem establishes the existence of a fake $3\CP^2\# 5\bCP^2$ which contains a square zero symplectic torus with simply connected complement.  In the statement $X$ denotes the manifold of Theorem \ref{cool} and $M$  the manifold   of Proposition \ref{SxT}. Recall that $X$ contains a square zero genus 2 surface $F$ with  simply connected complement. Also,  $M$ contains a square zero genus two surface, also denoted $F$,  and four homologically essential Lagrangian tori.

 \begin{thm}\label{10baby} One can do Luttinger surgeries on the symplectic  sum of $X$ with $M$  to obtain a minimal symplectic manifold $B$ homeomorphic but not diffeomorphic to $3\CP^2\# 5\bCP^2$.  This manifold contains a pair of tori with trivial normal bundles    $T_3, T_4$, with $B-(T_3\cup T_4)$ simply connected. The tori $T_3$ and $T_4$   can be assumed to be either Lagrangian or symplectic.

 \end{thm}

\begin{proof}  Call the four Lagrangian tori in the sum
$$A= X-nbd(F)\cup M-nbd(F)$$
$T_1, T_2, T_3, T_4$ (these lie in the $M$ side).  Since  $X-nbd(F)$ is simply connected, and the boundary of $M-nbd(F)$ carries the loops $a_1,b_1,a_2,b_2$, $\pi_1(A)$ is generated by the  loops $x$ and $y$ in $M$.  Proposition \ref{SxT} shows that the meridians of $T_1$  and $T_2$ are  trivial. Then $-1/1$ Luttinger surgery along $m_1$ in $T_1$ and $m_2$  in $T_2$ (and no surgery, i.e. $1/0$, surgery on $T_3$ and $T_4$) yields a simply connected, symplectic manifold $B$ with $e=10$ and $\sigma=-2$.

Notice that $\pi_1(B-(T_3\cup T_4))=1 $.  The Lagrangian torus $T_3$ is homologically essential in $B$, since it intersects the $y\times b_2$ torus transversally once.   Thus one can perturb the symplectic structure on $B$   so that $T_3$ becomes symplectic. Similarly for $T_4$.

 The manifold $X$ is minimal by Theorem \ref{cool}. Since it is simply connected, it is not an $S^2$ bundle over $F$.   The symplectic manifold $M'$ obtained from the Luttinger surgeries on $T_1$ and $T_2$ in $M$ is minimal by the same argument as in the proof of Theorem \ref{five}.
Since $A$ is the fiber sum of $X$ and $M'$,   Usher's theorem implies that $A$ is minimal.
\end{proof}

The nullhomologous tori  in $B$ coming from $T_1$ and $T_2$  can be used  to produce infinitely many non-diffeomorphic smooth manifolds homeomorphic to $3\CP^2\# 5\bCP^2$, using the  argument of Corollary \ref{infinitely}. Moreover, the tori $T_3$ and $T_4$ with simply connected complement  can be used to perform knot surgery in the sense of \cite{FS6}.

\bigskip

More  simply connected examples can be obtained by taking the symplectic sums with $V$ or $B$. For example, the symplectic  sum of $B$ with  $W_1$  along a parallel copy of $F$ and  $F_1$ is simply connected by Lemma \ref{ontolem}.  It is minimal, symplectic, and  homeomorphic but not diffeomorphic to $5\CP^2\# 11\bCP^2$. Similarly the symplectic sum of $B$ with $W_2$ gives a  minimal symplectic
 exotic  $5\CP^2\# 9\bCP^2$,  the symplectic sum of $V$ with $W_1$ gives a  minimal symplectic
 exotic
$3\CP^2\# 11\bCP^2$, and  the symplectic sum of $V$ with $W_2$ gives a  minimal symplectic
 exotic
$3\CP^2\# 9\bCP^2$ (c.f \cite{dpark}).

\medskip

This process can   be iterated. The following corollary gives one such example.

\begin{cor}\label{family} Given non-negative integers $m$ and $n$, let $C_{m,n}$ denote the
 symplectic sum of $X$ with $m$ copies of $W_1$ and $n$ copies of $W_2$ along parallel copies of $F$. Then $C_{m,n}$ is a simply connected   minimal symplectic manifold   homeomorphic but not diffeomorphic to $(1+2m+2n)\CP^2\#(3+6m+4n)\bCP^2$.
\end{cor}
 \begin{proof} Lemma  \ref{ontolem} and induction proves that $C_{m,n}$ is simply connected.  We have
 $$e(C_{m.n})= e(X)+me(W_1) + ne(W_2)+4(m+n)= 6+8m+6n$$
 and
 $$\sigma(C_{m,n})= \sigma(X) +m\sigma(W_1)+ n\sigma(W_2)= -2(1+2m+n).$$

The intersection form of $C_{m,n}$ is odd for any $n$. To see this, consider the class represented  by $\hat H\times \{s\}$ in $M=\hat H\times \Sigma$. We form $C_{m,n}$ by taking the symplectic  sum of $M$ with $m$ copies of $W_1$ and $n+1$ copies of $W_2$ along $m+n+1$ parallel copies of $F$, say $\{p_i\}\times\Sigma,  i=1,\cdots m+n+1$. We can arrange that the first symplectic sum (of $M$ with $W_2$) is a relative symplectic sum  (\cite{Gompf}) which lines up the punctured $\hat H$ with one of the exceptional spheres in $W_2$. This produces a torus with square $-1$ in $X$. We can then take relative symplectic sums so that this torus lines up with vertical (square zero) tori in $W_1=T^2\times S^2\#4\bCP^2$  and
$W_2=T^4\#2\bCP^2$ (rather than the exceptional curves), so that $C_{m,n}$ contains an embedded  genus $m+n+1$ surface of square $-1$. Hence The intersection form of $C_{m,n}$ is odd.

 By Freedman's theorem $C_{m,n}$  is
 homeomorphic  $(1+2n)\CP^2\#(3+4n)\bCP^2$.
Minimality is proved as before using Usher's theorem.    Taubes's results then imply that $C_{m,n}$ is not diffeomorphic to $(1+2m+2n)\CP^2\#(3+6m+4n)\bCP^2$.
 \end{proof}

There are many more  minimal symplectic manifolds that can be produced using Theorem \ref{prop1}.
One can fill out most of  the region $\{(c_1,\chi_h)\ | \ c_1^2<8\chi_h, c_1^2\text{ even }\}$ using these constructions.  We will address this in another article.

\medskip

The following theorem is useful to construct manifolds with different smooth structures.

\begin{thm}\label{B1} There exists a minimal symplectic manifold $B_1$ containing a pair of Lagrangian tori $T_3$ and $T_4$ and a square zero symplectic genus 2 surface $F$ so that $T_3,T_4$ and $F$ are pairwise disjoint and
\begin{enumerate}
\item $\pi_1(B_1-(F\cup T_3\cup T_4))= \ZZ^2$, generated by $t_1$ and $t_2$.
\item The inclusion  $B_1-(F\cup T_3\cup T_4)\subset B_1$ induces an isomorphism on fundamental groups. In particular the meridians $\mu_F,\mu_3,\mu_4$ all vanish in $\pi_1(B_1-(F\cup T_3\cup T_4))$.
\item The Lagrangian push offs  $m_3,\ell_3$ of $\pi_1(T_3)$ are sent to $1$ and $t_2$ respectively in the   fundamental group  of $ B_1-(F\cup T_3\cup T_4)$.
\item The Lagrangian push offs  $m_4,\ell_4$ of $\pi_1(T_4)$ are sent to $t_1$ and $t_2$ respectively in the   fundamental group  of $ B_1-(F\cup T_3\cup T_4)$.
\item The push off $F\subset   B_1-(F\cup T_3\cup T_4)$ induces a map on fundamental groups with image the subgroup generated by  $t_2$.
\end{enumerate}
The symplectic structure may be perturbed so that one or both of the tori $T_3$, $T_4$ are symplectic.
 \end{thm}
\begin{proof}     Construct  $B_1$   by  starting with the manifold $Z$ of Theorem \ref{prop1}, and doing Luttinger surgeries along four of the Lagrangian tori, just as was done in the proof of Theorem \ref{cool}. Explicitly, we do $1/1$ Luttinger surgery on $T_1'$ along $m_1'$, $-1/1$ surgery on $T_1$ along $m_1$, then  $-1/1$ surgery on $T_2$ along $\ell_2$, and finally $1/1$ surgery on $T_2'$ along $\ell_2'$.

The resulting manifold $B_1$ has two remaining Lagrangian tori, $T_3$ and $T_4$, and the same calculation as in the proof of Theorem \ref{cool} shows that $\pi_1(B_1-(F\cup T_3\cup T_4))$ is generated by $a_2$ and $y$ and $[a_2,y]=1$.  Moreover, the meridian $\mu_F=1$ and also the other two meridians $\mu_3=\mu_4=1$. The Lagrangian push offs are given by
$$m_3=1, \ell_3=a_2, m_4=y, \ell_4=a_2$$
using Theorem \ref{prop1}.

Thus $\pi_1(B_1-(F\cup T_3\cup T_4))$ is a quotient of $\ZZ a_2\oplus \ZZ y$. But Proposition \ref{sixtori}  shows that $H_1(B_1-(F\cup T_3\cup T_4))=\ZZ\oplus \ZZ$, and so $\pi_1(B_1-(F\cup T_3\cup T_4))=\ZZ a_2\oplus \ZZ y$.

The tori $T_3$ and $T_3$ have dual tori by Proposition \ref{sixtori},  and so they are homologically essential and linearly independent in second homology. Hence the symplectic form can be perturbed slightly so that one or both of $T_3$ and $T_4$ are symplectic (\cite{Gompf}).

Setting $t_1=y$ and $t_2=a_2$ completes the proof.
\end{proof}

The following corollary shows that at the expense of a small stabilization, any symplectic manifold containing a square zero symplectic genus 1 or 2 surface  can be used to produce infinite families of manifolds with different Seiberg-Witten invariants.   Note that if the surface $G$ has genus 1, then a symplectic sum with $E(1)$ and the knot surgery procedure of Fintushel-Stern \cite{FS6} gives a similar result, at the cost of $12$ to $e$ and $-8$ to $\sigma$.  Thus this result can be viewed as an extension (it includes genus 2) and an improvement (the cost to $e$ is 6 or 10). (See also \cite{FS4}.)

\begin{cor}\label{infinitely2} Suppose $M$ is any symplectic manifold which contains a square zero genus 1 or 2 symplectic surface  $G$ such that the inclusion $G\subset M$ induces the trivial map on fundamental groups.  Then there exists an infinite  family of smooth manifolds  $M_n$ so that
$\pi_1(M_n)=\pi_1(M)$,
 $e(M_n)= e(M)+ 2+ 4\text{ genus}(G)$,
 $\sigma(M_n)=\sigma(M)-2$, and the Seiberg-Witten invariants of  $M_n$ are different from the Seiberg-Witten invariants of $M_m$ if $m\ne m$.  Moreover $M_1$ is symplectic.
\end{cor}

\begin{proof} This is an application of the theorem of Fintushel and Stern \cite{FS5}, which is based on the gluing formula for Seiberg-Witten invariants of Mrowka-Morgan-Szabo \cite{MMS}, and is similar to the argument given in Corollary \ref{infinitely}  above.

Starting with $M$, take the symplectic sum of $M$ with the manifold $B_1$ of Theorem \ref{B1} along $F$ if $G$ has genus 2, and along $T_3$ (after making $T_3$ symplectic by perturbing the symplectic structure) if $G$ has genus 1. Call the result $S$. Then $S$ is a symplectic manifold satisfying $e(S)=e(M)+ 6 + 4(\text{genus}(G)-1)$ and  $\sigma(M_n)=\sigma(M)-2$.

The Seifert-Van Kampen theorem  and Theorem \ref{B1} implies that  $\pi_1(S)=\pi_1(M)*\ZZ t_1$. The inclusion $S-T_4\subset S$ induces an isomorphism on fundamental groups. The Lagrangian push   offs of $T_4$ are $m_4=t_1$ and $\ell_4=1$, and the meridian is $\mu_4=1$.

We follow the argument of Corollary \ref{infinitely}.  Let $M_0$ denote the manifold obtained from $S-nbd(T_4)$ by torus  filling  by $T^2\times D^2$ in such a way that $\al=S^1 \times \{(1,1)\}$ is sent to $\ell_4$, $\be=\{1\}\times S^1 \times \{1\}$ is sent to $\mu_4$, and $\mu=\{(1,1)\}\times \partial D^2$ is sent to $m_4^{-1}$.   Then $\pi_1(M_0)=\pi_1(M)$ since $t_1$ is killed  This is not a Luttinger surgery, so $M_0$ need not be
symplectic.

Let $\tilde{T}\subset Y_0$ denote the resulting core torus. Then $Y_0-nbd(T)=S-nbd(T_4)$.
With the coordinates $\al, \be, \mu$, $S$ is obtained from $M_0$ by $0/1$ surgery on $\tilde T$ along $\be=\mu_4$. Moreover, $\be$ is nullhomologous in $M_0-\tilde T$ since $\mu_4=0$ in $H_1(B_1-(F\cup T_3\cup T_4))$.  Note that $S$ is a symplectic manifold with $b^+>1$, since $T_3$ has a dual torus (namely
$b_2\times y$) and either $T_4$ also has a dual torus $x\times b_2$  (in case $G$ has genus 2) or $F$ has a dual torus $\hat H\times \{p\}$.  Hence the Seiberg-Witten invariants of $S$ are non-zero.

Since the surgery on $T_4$ producing $M_0$ from $S$  kills the generator $t_i$ of $H_1$, it also kills the hyperbolic pair containing $T_4$, i.e. $\tilde T=0$ in $H_2(M_0)$.

Let $M_n$ denote the manifold obtained from $1/n$ surgery on $\tilde T$ in $M_0$ along $\be$. Note that $M_1$ can also be viewed as $-1/1$ Luttinger surgery on $T_4$ in $S$ along $m_4$.  Thus $M_1$ is symplectic.

The main result of \cite{FS5} then shows that the family $M_n $ obtained from  $1/n$ surgery  on $\tilde T$ along $\be$ contains infinitely many diffeomorphism types.

\end{proof}

 \section{Non-trivial fundamental group}

\subsection{Fundamental group $\ZZ$}

We turn now to a useful example of a symplectic 4-manifold $X_1$ with fundamental group $\ZZ$.

\begin{thm} \label{ZZ} There exists a  minimal symplectic 4-manifold $X_1$   with $\pi_1(X_1)=\ZZ$, $e(X_1)=6$ and $\sigma(X_1)=-2$. Moreover, $X_1$ contains a symplectic surface $F$ of genus 2  and square zero and symplectic  torus $T$ of square zero disjoint from $F$ so that $\pi_1(X_1-(T\cup F))\cong \ZZ$.  The homomorphisms $\pi_1(X_1-(T\cup F))\to \pi_1(X_1-F)$  and $\pi_1(X_1-T)\to \pi_1(X_1)$ induced by inclusion are isomorphisms.

 Furthermore, the homomorphism $\pi_1(T)\to \pi_1(X_1-F)$ induced by inclusion takes one generator $ t_1\in \pi_1T$ to the generator of $\pi_1(X_1-F)$ and the other $t_2$ to the trivial element $1$. The homomorphism $\pi_1(F)\to \pi_1(X_1)$ induced by inclusion is trivial.
\end{thm}
\begin{proof}  We follow the beginning of the  proof of Theorem   \ref{cool},   performing   Luttinger surgeries on the  Lagrangian tori in $Z-F$, where $Z$ is the manifold of Equation (\ref{Zmanifold}).  We refer to the calculations of Theorem \ref{prop1}.

First $1/1$ surgery on $T_1'$ along $m'_1$ imposes the relation $b_1=[a_2^{-1}, a_1^{-1}]$. Then $-1/1$ surgery on $T_1$ along $m_1$ imposes the relation $x=[b_1^{-1}, y^{-1}]= [[a_1^{-1}, a_2^{-1}], y^{-1}]=1$. The last equality comes from Equation
(\ref{eq4.2}).

Next perform $-1/1$ Luttinger surgery on $T_2$ along $\ell_2$. This imposes
$a_1=[x^{-1}, b_1]=1$ and hence also $b_1=1$.  Then $1/1$ surgery on $T_2'$ along $\ell_2'$ sets
$b_2=[b_1,a_2]=1$.   A $-1/1$ surgery on $T_3$ along $ \ell_3$
sets $a_2=[b_2^{-1},y^{-1}]=1.$

This leaves $T_4$ untouched.
Call the resulting manifold $X_1$. Thus the classes $x,a_1,b_1,a_2,$ and $b_2$ equal $1$ in $\pi_1(X_1-(T_4\cup F))$. Since $H_1(Z)=\ZZ^6$ and we performed five Luttinger surgeries, $\pi_1(X_1-(T_4\cup F))$ is generated by $y $ and $H_1(X_1)=\ZZ$, so that $\pi_1(X_1-(T_4\cup F))=\pi_1(X_1-F)=\pi_1(X_1)=\ZZ y$.

Then the torus $T_4$ carries the classes $y$ and $a_2$, and the surface $F$ is generated by the classes $a_1,b_1,a_2,b_2$, and  so the assertions about fundamental groups follow.  Note that $e(X_1)=e(Z)=6$ and $\sigma(X_1)=\sigma(Z)=-2$.

Lastly, the torus $T_4$ meets the torus  $Q=x\times b_2$ in $M$ transversally in one point. Since $Q$ misses all other tori  in the construction and $F$, it survives to provide $T_4$ with a dual class in $H_2(X_1)$. In particular, $[T_4]\ne 0$ in $H_2(X_1)$.  Hence  the symplectic form on $X_1$ can be perturbed slightly so that $T_4$  (which we rename $T$) becomes symplectic (\cite{Gompf}). We relabel its generators $t_1=y$ and $t_2=a_2$.

Minimality follows just as in the proof of Theorem \ref{five}.
\end{proof}

\medskip

The intersection form of $X_1$ is equivalent to $2(1)\oplus 4(-1)$. But $X_1$ is not diffeomorphic to $(S^1\times S^3)\#2\CP^2\#4\bCP^2$,  since   $X_1$ is a minimal symplectic manifold with $b^+=2>1$, and hence has non-vanishing Seiberg-Witten invariants \cite{taubes3}.  On the other hand the Seiberg-Witten invariants of the connected sum $(S^1\times S^3)\#2\CP^2\#4\bCP^2$ must vanish. It is an interesting question whether $X_1$ is homeomorphic to $(S^1\times S^3)\#2\CP^2\#4\bCP^2$.

\medskip

Our  interest in the manifold $X_1$ is two fold. First, it is the smallest known (to us) symplectic 4-manifold with fundamental group $\ZZ$, where  we measure the size using the Euler characteristic (or equivalently the second Betti number).  (Constructions of symplectic manifolds with fundamental group $\ZZ$ can be found in the literature, e.g. \cite{OS}, \cite{Gompf}, \cite{smith}.) The other reason is that it can be used as a smaller replacement  for the elliptic surface $E(1)$ typically used to control fundamental groups of symplectic 4-manifolds. We will illustrate this in   the following theorem, which also refers to the manifold $B$ constructed in Theorem \ref{10baby}.

\begin{thm} \label{thm4} Let $L$ be a symplectic 4-manifold containing a symplectic torus $T'$ with trivial normal bundle such that $\al,\be\in \pi_1(L)$ represent   the two generators of $\pi_1(T')$. Then

\begin{enumerate}

\item The symplectic sum of $X_1$ and $L$ along $T$ and $T'$, $X_1\#_T L$, using the gluing map $\phi:T\to T'$
inducing $$t_1\mapsto \al, t_2\mapsto \be$$
admits a symplectic structure which agrees with that of $X_1$ and $L$ away from $T,T'$  and satisfies
$$e(X_1\#_T L)=e(L)+6,\ \sigma(X_1\#_T L)=\sigma(L)-2,\ \text{ and } \ \pi_1(X_1\#_T L)=\pi_1(L)/N(\be)$$
where $N(\be)$ denotes the normal subgroup of $\pi_1(L)$ generated by $\be$.

 \item  The symplectic sum of $B$ and $L$ along $T_3$ and $T'$, $B\#_T L$ admits a symplectic structure which agrees with that of $B$ and $L$ away from $T_3,T'$  and satisfies
$$e(B\#_T L)=e(L)+10,\ \sigma(B\#_T L)=\sigma(L)-2,\ \text{ and } \ \pi_1(B\#_T L)=\pi_1(L)/N(\al,\be).$$

\item Suppose $L$ is a symplectic 4-manifold and $\tilde{F}\subset L$ is a genus 2 symplectic surface with trivial normal bundle. Then the symplectic sum $L\#_F X$ of $L$ and $X$ along their genus 2 surfaces is a symplectic manifold satisfying
$$e(L\#_F X)=e(L)+10, \ \sigma(L\#_F X)=\sigma(L) -2$$
and $$\pi_1(L\#_F X)=\pi_1(L)/N$$
where $N$ denotes the normal subgroup of $\pi_1(L)$ generated by the image of $\pi_1(\tilde{F})\to \pi_1(L)$.
\end{enumerate}

\end{thm}

\begin{proof} The
assertions about $e$ and $\sigma$ are straightforward.
 The fundamental group assertion is proved using the Seifert-Van Kampen theorem.
\end{proof}

Theorem \ref{thm4} can be restated informally by saying  that at a cost of $6$ to the Euler characteristic one can symplectically kill one class in the fundamental group of a symplectic manifold, provided that class is carried by a symplectic torus.  Similarly at a cost of $10$ to $e$ one can kill two classes carried on a symplectic torus.  Lastly   at a cost of $10$ to $e$ one can kill four classes  carried by a symplectic genus 2 surface.

\subsection{Arbitrary fundamental group}
The fundamental group $G$ of a closed, orientable 4-manifold $M$ determines all its Betti numbers except $b_2$. Moreover, $b_2(M)\ge b_2(G)$, and hence
\begin{equation}\label{lowerbound}
e(M)= 2-2b_1(M)+b_2(M)= 2-2b_1(G)+ b_2(M)\ge 2-2b_1(G)+b_2(G)
\end{equation}

One can get more subtle lower bounds on $e(M)$ by studying the algebraic topology (e.g.~ the ring structure) of $K(G,1)$; see   \cite{KL1,KL2}.   In particular, thinking of the rank of $H_2$ as a measure of the size of a 4-manifold,  one sees that minimizing the Euler characteristic is the same as minimizing this  size, among manifolds with a given fundamental group.

\bigskip

As explained in Section 4 of \cite{BK},   the existence of the symplectic manifold $X_1$ and its symplectic torus $T$ of Theorem \ref{thm4} allows us to   improve  (by 50\%) the main result of \cite{BK} to the following theorem.

\begin{thm}\label{50percent} Let  $G$ be a finitely presented group that has a   presentation with $g$
generators and $r$ relations. Then
   there exists a symplectic 4-manifold $M$ with $\pi_1M\cong G$,  Euler
characteristic
$e(M)=10+ 6(g+r)$, and signature $\sigma(M)=-2-2(g+r)$.   \qed
   \end{thm}
 \begin{proof}
The proof relies on our construction in \cite[Theorem 6]{BK} of  a symplectic 4-manifold $N$ is constructed whose fundamental group  contains   classes $s,t, \gamma_1,\cdots,\gamma_{r+g} $  so that $$G\cong \pi_1(N)/N(s,t, \gamma_1,\cdots,\gamma_{r+g})$$  where
$N(s,t, \gamma_1,\cdots,\gamma_{r+g}) $  denotes the normal subgroup generated by the  classes $s,t, \gamma_1,\cdots,\gamma_{r+g} $.

  Moreover, $N$ contains symplectic tori $T_0$, $T_1,\cdots, T_{g+r}$ so that the two generators of $\pi_1(T_0)$ represent $s$ and $t$,  and for $i\ge 1$ the two generators of $\pi_1(T_i)$ represent $s$ and $\gamma_i$.  The manifold $N$ satisfies $e(N)=0$ and $\sigma(N)=0$; in fact $N$ is a product $Y\times S^1$ where $Y$ is a 3-manifold that fibers over $S^1$.

 Let $B$ denote the manifold of Theorem \ref{10baby}.
Take the fiber sum of $N$ with $B$ along $T_0$, and  $g+r$ copies of the manifold $X_1$ of Theorem \ref{ZZ} along the tori $T_i$, $i\ge 1$  using an appropriate gluing map as in Theorem \ref{thm4},
so that $s,t$ and the $\gamma_i$ are killed.  Then a repeated application of Theorem \ref{thm4} computes
$$\pi_1(M)=\pi_1(N)/N(s,t,\gamma_i)=G,\ e(M)=10 +6(g+r), \ \sigma(M)=-2-2(g+r).$$
\end{proof}

 An examination of the proof of Theorem 6 of \cite{BK} shows that Theorem \ref{50percent} can be improved for certain presentations, namely, one can find $M$ so that $e(M)=10 + 6(g'+r)$ and $\sigma(M)=-2-2(g'+r)$, where $g'$ is the number of generators which appear in some relation with negative exponent. Thus if $G$ has a presentation with $r$ relations in which every generator appears only with positive exponent in each relation, then there exists  a symplectic $M$ with $\pi_1(M)=G$, and $e(M)=10+6r,\ \sigma(M)=-2-2r$.  Moreover, using Usher's theorem \cite{usher} one sees that
the manifolds constructed are minimal.

For example, if $G$ is the free product of $n$ finite cyclic groups, then $G$ has a presentation
$$G=\langle x_1,\cdots, x_n\ | \ x_1^{a_1}, \cdots, x_n^{a_n}\rangle$$
with all the $a_i>0$, and so there exists a symplectic 4-manifold with fundamental group $G$ and $e=10+6n$.

\medskip

If one uses $E(1)$ instead of $B$ in the proof of Theorem \ref{50percent}, the resulting manifold has $e=12+6(g+r)$ and $\sigma=-8-2(g+r)$, and contains   a symplectic torus which lies in a cusp neighborhood. Thus the geography results of J. Park \cite{park} can be improved to find a larger region of the $(c_1^2,\chi_h)$ plane for which to each pair of integers in that region one can find   infinitely many non-diffeomorphic, homeomorphic minimal symplectic manifolds with fundamental group $G$.

 \smallskip

 Another (decidedly minor) improvement concerns groups of the form $G\times\ZZ$: for a presentation of $G$ as above there exists a symplectic manifold $M$ with $\pi_1(M)=G\times\ZZ$, $e(M)=  6(g'+r+1)$, and $\sigma(M)=-2(g'+r+1)$. The reason is that one step  in the proof of \cite[Theorem 6]{BK} consists of taking a symplectic sum with $E(1)$ to kill two generators ($t$ and $s$ along $T_0$ in the notation of the proof of Theorem \ref{50percent}). But to get $G\times \ZZ$ it suffices to kill $t$, for which the manifold $X_1$   can be used instead of $B$, using Theorem \ref{thm4}. This only adds $6$ to $e$.

 Notice that if $G$ is the fundamental group of a 3-manifold $Y$ that fibers over $S^1$, then $Y\times S^1$ is a symplectic 4-manifold with fundamental group $G\times \ZZ$ and Euler characteristic zero.

\bigskip

Suppose $M$ is a symplectic 4-manifold containing  a symplectic torus $T$ with trivial normal bundle so that $\pi_1(M-T)=1$ or $\ZZ$ and so that pushing $T$ into $M-T$ induces a surjection $\pi_1(T)\to \pi_1(M-T)$.  One can prove, using the adjunction formula \cite{OS2},  that such an $M$ must have $b^+>1$ and $b^->0$, and hence if $M$ is simply connected $e(M)\ge 6$.
The manifold $B$ of Theorem \ref{10baby}  is such a simply connected manifold and has $e(B)=10$.
 If $\pi_1(M)=\ZZ$ such an $M$ must have $e(M)\ge 3$.   The manifold $X_1$ has $e(X_1)=6$.

 Further improvements in the geography problem for symplectic manifolds will be obtained if such an $M$ is found with $6\leq e(M)< 10$ in the simply connected case and $3\leq e(M)<6$ in the $\ZZ$ case.  The search for such a manifold is a promising direction for future study.

\subsection{Free groups}

\begin{thm}\label{free} For any $n$ there exists a symplectic 4-manifold $D_n$ with fundamental group free of rank $n$, $e(D_n)=10$, and $\sigma(D_n)=-2$.
\end{thm}
\begin{proof} This is explained in \cite[Theorems 8 and 11]{BK}, but here,  instead of taking a symplectic sum with $E(1)$, one uses the smaller manifold $B$ of Theorem \ref{10baby}.  Specifically, let $E$ be a closed surface of genus $n$, with standard generators $x_1,y_1,\cdots, x_n,y_n$ of $\pi_1(E)$. Let $D:E\to E$ be the diffeomorphism given by composite of the $n$ Dehn twists along the $x_i$ and let $Y$ be the corresponding 3-manifold which fibers over $S^1$. Then $Y\times S^1$ is a symplectic 4-manifold containing a symplectic torus $T'=S\times S^1$ where $S$ is a section of $Y\to S^1$.

The standard calculation of $\pi_1(Y)$ as an HNN extension shows
$$\pi_1(Y\times S^1)=\langle x_i,y_i, t\ | \ t x_i t^{-1}=x_i, t y_i t^{-1}= y_ix_i\rangle\oplus \ZZ s.$$

By Theorem \ref{thm4},  taking the symplectic sum of $Y\times S^1$ with $B$ along $T'$ and $T$  has fundamental group obtained by killing $t$ and $s$ in $\pi_1(Y\times S^1)$, which also kills the $x_i$ (since the relation $y_i=y_ix_i$ implies that $x_i=1$), leaving the free group  generated by the $y_i$. Theorem \ref{thm4} together with  $e(Y\times S^1)=0$ and $\sigma(Y\times S^1)=0$ gives the result.
\end{proof}

Notice that the smooth 4-manifold obtained by taking connected sums of $n$ copies of $S^1\times S^3$  has fundamental group free of rank $n$ and $e=2(1-n)$.  Kotschick has shown (\cite{kot2}) that any symplectic 4-manifold with fundamental group free of rank $n$ has $e\ge  \frac{6}{5}(1-n)$. Thus the gap in size (measured say by the rank of the second homology) between the smallest smooth and symplectic manifolds must grow  linearly with the rank.     The manifold of Theorem \ref{free}, with fundamental group free of rank $n$ and $e=10$, is the smallest known symplectic 4-manifold with free fundamental group. With the exception of rank 1 (for which the manifold $X_1$ has $e=6$)  we do not know of any symplectic 4-manifold $M$ with fundamental group free of rank $n$ which satisfies
$$\tfrac{6}{5}(1-n)\leq e(M)< 10.$$

\bigskip

\subsection{Fundamental groups associated to   surface bundles over the circle}

By exactly the same proof as Theorem \ref{free} one establishes (see \cite[Theorem 8]{BK}).

\begin{thm} If $H:F\to F$ be an orientation preserving diffeomorphism of a closed orientable surface, fixing a point $z$, and $G$ is the quotient of $\pi_1(F,z)$ by the normal subgroup generated by the words $x^{-1}H_*(x)$, then there exists a symplectic 4-manifold with fundamental group $G$, Euler characteristic $e=10$, and signature $\sigma=-2$.

Moreover, there exists a symplectic 4-manifold with fundamental group $G\times \ZZ$, Euler characteristic $e=6$, and signature $\sigma=-2$.
\qed
\end{thm}
We can also produce small symplectic 4-manifolds with the same fundamental group as a fibered 3-manifold.

\begin{thm} If $H:F\to F$ be an orientation preserving diffeomorphism of a closed orientable surface and $Y\to S^1$ the $F$-bundle over $S^1$ with monodromy $H$, then there exists a symplectic 4-manifold
$D$ with $\pi_1(D)=\pi_1(Y)$, $e(D)=6$, and $\sigma(D)=-2$.
\end{thm}

\begin{proof} The manifold $Y\times S^1$ admits a symplectic structure so that the torus $S\times S^1$ is symplectic, where $S\subset Y $ is a section.   Theorem \ref{thm4} shows that taking a symplectic sum of $Y\times S^1$ with $X_1$ along an appropriate diffeomorphism of tori yields a symplectic manifold $D$ in which the homotopy class of the $S^1$ factor is killed, and so $\pi_1(D)=\pi_1(Y)$. This manifold has $e=6$ and $\sigma=-2$.
\end{proof}

In \cite{KL2} it is established that for {\em any} closed 3-manifold group $G$  and any number $\sigma$, there exists a smooth 4-manifold with fundamental group $G$ and  $e=2+|\sigma|$, and that this is the smallest possible Euler characteristic among all 4-manifolds with fundamental group $G$ and signature $\sigma$.   In particular,   there is a smooth 4-manifold with fundamental group $G$,   $\sigma=-2$, and $e=4$.  We have found a symplectic 4-manifold whose Euler characteristic is within $2$ of the smooth minimum for those $G$ which are the fundamental group of a fibered 3-manifold.

\subsection{Free abelian groups}

The manifold $X_1$ constructed in Theorem \ref{ZZ} is currently the smallest known  symplectic 4-manifold with fundamental group infinite cyclic, with $e=6$.
Producing small smooth 4-manifolds with free abelian fundamental groups poses an interesting challenge \cite{KL1,KL2}.  The number of required relations in a presentation grows quadratically in the number of relations, and so one expects many 2-handles in a handlebody presentation.

Finding symplectic examples is harder. For free abelian groups of even rank  a nearly complete answer was found in the collection of symplectic (in fact K\"ahler) manifolds $Sym^2(F_n)$ ($F_n$ a surface of genus $n$). The manifold $Sym^2(F_n)$   has fundamental group $\ZZ^{2n}$ and  minimizes the Euler characteristic  among symplectic manifolds with fundamental group $\ZZ^{2n}$ except possibly when $n\equiv 2\pmod{4}$ (\cite{BK}).

 For odd rank free abelian fundamental group the situation is less clear. We do not know if there exists a symplectic manifold $M$ with fundamental group $\ZZ$  and
 $3\leq e(M)<6$, nor any reason why such a manifold cannot exist.
 For rank 3, we have the following result.

 \begin{cor}\label{Z3} There exists a symplectic 4-manifold $X_3$ with fundamental group $\ZZ^3$, $e=6$ and $\sigma=-2$.
  \end{cor}
  \begin{proof} Start with the manifold $Z$ of Theorem \ref{prop1}. We do a sequence of Luttinger surgeries, as in the proof of Theorems \ref{cool}  and \ref{ZZ}.  We first do $1/1$ surgery on $T_1'$ along $m_1'$ to impose the relation $b_1=[a_2^{-1},a_1^{-1}]$. Then $-1/1$ surgery on $T_1$ along $\ell_1$ yields $a_1=[b_1^{-1}, y^{-1}]=[[a_2^{-1},a_1^{-1}], y^{-1}]=1$ using Equation (\ref{eq4.2}). Thus also $b_1=1$. Next, $1/1$ surgery on $T_2'$ along $\ell_2'$ gives $b_2=[b_1,a_2]=1$. The remaining three generators $x,y, a_1$ commute using Equation (\ref{eq4.2}) and the fact that $[x,y]=1$ (We do not remove $F$). Call the result $X_3$.

  Since each surgery decreases the first Betti number by one, and $H_1(Z)=\ZZ^6$  it follows that $\pi_1(X_3)=H_1(X_3)=\ZZ^3$.
\end{proof}

 The smallest previously known example of a symplectic 4-manifold with fundamental group $\ZZ^3$ has $e=12$. Any such symplectic manifold must have $e\ge 3$ (\cite{BK}).   We know no reason why one cannot exist.

 \medskip

 More generally, the technique of  \cite[Theorem 20]{BK}  allows us to improve  the construction  of symplectic 4-manifolds with odd rank free abelian groups by taking the fiber sum of $Sym^2(F_n)$  with the manifold $X_1$ along an appropriate torus, rather than the larger manifold $K$ of \cite[Lemma 18]{BK}. We refer the interested reader to  \cite{BK} for details of the proof of the following corollary, which follows by replacing every occurrence of the symbol $K$ by $X_1$ in the proof of \cite[Lemma 18]{BK}.

 \begin{cor}\label{odd} There exists a symplectic 4-manifold $M$ with $\pi_1(M)=\ZZ^{2n-1}$ such that $
 e(M)=9-5n +2n^2$ and $\sigma(M)=-1-n$.\qed
\end{cor}

This is an improvement over the best previously known example,($e(M)= 15-5n +2n^2$ of
\cite[Theorem 20]{BK}), but still far from the best-known lower bound $6-7n+2n^2\leq \min_{\pi_1(M)=\ZZ^{2n-1}}e(M)$.

 \subsection{Abelian groups}
 Modifying the argument of Corollary \ref{Z3} gives the following.
  \begin{cor}\label{abelian} Given any $p,q,r\in \ZZ$ there exists a symplectic 4-manifold $X_{p,q,r}  $ with fundamental group $\ZZ/p\oplus\ZZ/q\oplus\ZZ/r$, $e=6$ and $\sigma=-2$.
  \end{cor}
  \begin{proof} Start with the manifold $X_3$ of Corollary \ref{Z3}.  Thus    $\pi_1(X_3)=\ZZ^3$, generated by $x,y$, and $a_1$, and $b_1=a_2=b_2=1$.  The Lagrangian tori $T_2,T_3$ and $T_4$ were not used to construct $X_3$. Then $1/p$ Luttinger surgery on $T_2$ along $\ell_2=a_1$ sets $a_1^p=1$. Similarly $1/q$ surgery on $T_3$ along $m_3=x$ and $1/r$ surgery on $T_4$ along $m_4=y$ sets $x^q=1$ and $y^r=1$.
\end{proof}

Note that the smallest previously known symplectic 4-manifolds with finite cyclic abelian group are  certain complex algebraic surfaces of general type with $e=10$ \cite{BPV}.   For sums of two or three abelian group the smallest previously known examples had $e=12$ \cite{Gompf, BK}.

\bigskip

The finitely generated abelian group  with $n$ generators
$$G=\ZZ^{n-k}\oplus \ZZ/d_1\oplus \ZZ/d_2\oplus \cdots\oplus \ZZ/d_k$$
(with $d_i>1$) has a presentation with $n$ generators, and $ {{n}\choose{2}}+k$ relations. Thus
Theorem \ref{50percent}  implies that there exists a symplectic manifold $N$ with fundamental group $G$ satisfying $$e(N)=10+6(n+ {{n}\choose{2}}+k), \sigma(N)=-2(n+ {{n}\choose{2}}+k +1)$$
The leading term of this expression for $e$ as a function of $n$ is  $3n^2$.  In other words, if we let
$$p(n)=\min\{e(N) \ | \ \pi_1(N) \text{ is abelian and is generated by $n$ elements}\}$$

one can say that
$$\lim_{n\to \infty}\frac{p(n)}{n^2}\leq 3.$$

Theorem 2 of \cite{KL1} implies that

$$\frac{1}{2}\leq  \lim_{n\to \infty}\frac{p(n)}{n^2}.$$

 The manifolds  $Sym^2(F_n)$ and those constructed in Corollary \ref{odd}   show that for  {\em free} abelian groups with $n$ generators, the leading term for $e$ is $\frac{n^2}{2}$. In fact, the following theorem shows this to be true for all finitely generated abelian groups.

 \begin{thm}\label{genabelian}  Assume $n$ is even by setting $d_k=1$ if necessary. There exists a symplectic manifold $N_G$ with fundamental group
 $$G=\ZZ^{n-k}\oplus \ZZ/d_1\oplus \ZZ/d_2\oplus \cdots\oplus \ZZ/d_k$$
 satisfying
 $$e(N_G)= \frac{1}{2} n^2 +\frac{19}{2}n +36 $$
 and
 $$ \sigma(N_G)=  -\frac{5}{2}n-8.$$
 \end{thm}
\begin{proof} Write $G= \ZZ/d_1\oplus \ZZ/d_2\oplus \cdots\oplus \ZZ/d_n$, where we allow $d_n$ to be any integer. By allowing $d_n=1$ if necessary, we  may assume that $n$ is even, say $n=2g-6$ for some $g\ge 3$.

Let $F_g$ denote an compact Riemann surface of genus $g$. The K\"ahler surface  (hence symplectic  4-manifold)
 $$S_g=Sym^2(F_g)=F_g\times F_g/\sim, \ (x,y)\sim(y,x)$$ has fundamental group isomorphic to $\ZZ^{2g}$.

A  pair $x,y$ of embedded curves in $F_g$ determine  a torus  $\tilde T(x,y)=x\times y$ in $F_g\times F_g$. This torus descends to an embedded torus   $T(x,y)\subset S_g$ when $x\cap y=\phi$. Moreover, in this case
$\tilde T(x,y)$ is Lagrangian in $F_g\times F_g$ and its image $T(x,y)\subset S_g$ is Lagrangian  (\cite[Proposition 21]{BK},  \cite{perutz}).

 Choose embedded curves $a_1,b_1,a_2,\cdots , a_g,b_g\subset F_g$   which represent a standard symplectic basis for   $H_1(F_g)$.  The composite
 $$F_g\cong F_g\times \{p\}\to F_g\times F_g\to S_g$$
 takes   $\{a_i,b_i\}$ to a basis for $H_1(S_g)$.
Let $x_4, y_4, x_5, y_5, \cdots, x_g, y_g$ denote $2g-6$ parallel copies of the curve $a_1$ on $F_g$. Consider the Lagrangian tori
\begin{equation}\label{eq6.1}T(a_1, b_2),  T(b_1, b_3),  T(a_2, a_3)\end{equation}
 and also
\begin{equation}\label{eq6.2}T(x_4, a_4), T(y_4,b_4), T(x_5,a_5),T(y_5,b_5),\cdots, T(x_g,a_g),T(y_g,b_g).\end{equation}

 These $2g-3$ tori are pairwise disjointly embedded in $S_g$.  They are  each homologically essential since each one intersects a dual torus of the same form (e.g.~$T(x_i,a_i)$ intersects  $T(b_1,b_i)$ transversally once).  Moreover, they are linearly independent since one can check that together with their  dual tori they span a hyperbolic subspace of $H_2(S_g)$.   Thus
 \cite[Lemma1.6]{Gompf} implies that   the symplectic form on $S_g$ can be perturbed slightly so that   these  $2g-3$ tori are symplectic.

 Since $\pi_1(S_g)$ is abelian, the Hurewicz map to $H_1(S_g)$ is an isomorphism, and  taking symplectic sums with $X_1$ has the same effect on the fundamental group as on the first homology. In particular, we need not worry about base points: killing a class in $H_1(S_g)$ also kills the corresponding class in $\pi_1(S_g)$.

 Now take the symplectic sum of $S_g$ with 3 copies of the manifold $B$ and $2g-6$ copies of the manifold $X_1$ as follows. The three copies of $B$ are summed along the three tori of Equation (\ref{eq6.1}). This kills $a_1,b_2, b_1, b_3, a_2$ and $a_3$.

 The $2g-6$ copies of $X_1$ are summed along the other tori  using appropriate diffeomorphisms. To start,   take the sum with $X_1$ along $T\subset X_1$ the diffeomorphism $T\to T(x_4,a_4)$ which takes the generators $t_1,t_2$ of $\pi_1(T)$ to $a_4^{-1},  x_4a_4^{d_1}$. Theorem \ref{thm4}  implies that this kills $x_4a_4^{d_1}$. Recall that $x_4$ is a parallel copy of $a_1$.  Since $a_1$ was killed before, this sets $a_4^{d_1}=1$. Continue by setting $b_4^{d_2}=1$, etc. Note that $d_i$ can be any integer, including $0$ or $1$.

This produces a symplectic manifold $N_G$ with fundamental group $G$.  Theorem \ref{thm4} shows that
 $$e(N_G)=30 + 6(2g-6)+ e(S_g)=12g-6 + (2g^2-5g+3)=\frac{1}{2} n^2 +\frac{19}{2}n +36 $$
 and
 $$ \sigma(N_G)= -2(2g-3) + \sigma(S_g)= -4g+6 +(1-g)=-\frac{5}{2}n-8.$$
\end{proof}

%**************************************************

%******* End of document **************************

%**************************************************


\bigskip

\bigskip

\begin{thebibliography}{99999}

\bibitem{ADK} D. Auroux, S.K. Donaldson, and L. Katzarkov, {\em Luttinger surgery along Lagrangian tori and non-isotopy for singular symplectic plane curves}. Math. Ann. 326 (2003), no. 1, 185--203.

\bibitem{A} A. Akhmedov, {\em Exotic smooth structures on $3\CP^2\#7\bCP^2$}. preprint (2006)
math.GT/0612130


\bibitem{AP} A. Akhmedov, B. Doug Park, {\em Exotic Smooth Structures on Small 4-Manifolds.} preprint (2007).

 \bibitem{B1} S. Baldridge, {\em New symplectic 4-manifolds with $b_+=1$.} Mathematische Annalen 333 (2005) 633-643.

 \bibitem{bald2} S. Baldridge, {\em An exotic symplectic structure on $\CP^2\#3\bCP^2$ and on $3\CP^2\#5\bCP^2$}. Preprint (2007).

\bibitem{BK} S. Baldridge and P. Kirk, {\em On symplectic 4-manifolds with prescribed fundamental group.} To appear in
Commentarii Mathematici Helvetici.

\bibitem{BK2} S. Baldridge and P. Kirk, {\em A symplectic manifold homemorphic but not diffeomorphic to
$\CP^2\#3\bCP^2$.} preprint math.GT/0702211.

\bibitem{BK3} S. Baldridge and P. Kirk, {\em Symplectic 4-manifolds with arbitrary fundamental group near the Bogomolov-Miyaoka-Yau line.} J. Symplectic Geom. 4 (2006), no. 1, 63--70.

\bibitem{BPV}  W. Barth,C. Peters, and A. Van de Ven,  ``Compact complex surfaces.'' Ergebnisse der Mathematik und ihrer Grenzgebiete (3), 4. Springer-Verlag, Berlin, 1984. x+304 pp.

%\bibitem{D} S. Donaldson, {\em
%Irrationality and the $h$-cobordism conjecture. }
%J. Differential Geom. 26 (1987), no. 1, 141--168.

% \bibitem{FS} R. Fintushel and R. Stern, {\em
%Invariants for Lagrangian tori.}
%Geom. Topol. 8 (2004), 947--968

%\bibitem{FS3} R. Fintushel and R. Stern, {\em Double node neighborhoods and families of simply %
% connected 4-manifolds with $ b^+=1$.}  J. Amer. Math. Soc. 19 (2006), no. 1, 171--180.

 \bibitem{FS2} R. Fintushel and R. Stern, {\em
Surgery on nullhomolgous tori and simply connected 4-manifolds with $b^+=1$.}
Preprint (2007) math.GT/0701044.

 \bibitem{FS5} R. Fintushel and R. Stern, {\em Seiberg-Witten invariants and infinite families of 4-manifolds.} preprint (2007).

  \bibitem{FS6} R. Fintushel and R. Stern, {\em   Knots, links, and $4$-manifolds.}
Invent. Math. 134 (1998), no. 2, 363--400.

 \bibitem{FS4} R. Fintushel and R. Stern, {\em
Families of simply connected 4-manifolds with the same Seiberg-Witten invariants.}
Topology 43 (2004), no. 6, 1449--1467.

\bibitem{FS7} R. Fintushel and R. Stern, {\em Will we ever classify simply-connected smooth 4-manifolds?.} Clay Mathematics Proceedings 5 (2006),  Floer Homology, Gauge Theory and Low Dimensional Topology, CMI/AMS Book Series, 225-240.

\bibitem{Freedman} M. Freedman, {\em The topology of four-dimensional manifolds.}
J. Differential Geom. 17 (1982), no. 3, 357--453.

\bibitem{Gompf} R. Gompf, {\em A new construction of symplectic manifolds}.  Ann. of Math. (2) 142 (1995), no. 3, 527--595.

\bibitem{GS} R. Gompf and A. Stipsitz, {\em 4-manifolds and  Kirby calculus}. Graduate Studies in Mathematics, 20. American Mathematical Society, Providence, RI, 1999.

\bibitem{KL1} P. Kirk and C. Livingston, {\em The Hausmann-Weinberger 4-manifold invariant of abelian groups.} Proc. Amer. Math. Soc. 133 (2005), no. 5, 1537--1546.

\bibitem{KL2} P. Kirk and C. Livingston, {\em The geography problem for 4-manifolds with specified fundamental group.} preprint math.GT/0608103.

\bibitem{Kot} D. Kotschick, {\em  On manifolds homeomorphic to $\CP^2 \#8\bCP^2$.}  Invent. Math. 95 (1989), no. 3, 591--600.

\bibitem{kot2} D. Kotschick, {\em Minimizing Euler characteristics of symplectic four-manifolds.}
 Proc. Amer. Math. Soc. 134 (2006), no. 10, 3081--3083.

\bibitem{Lut} K. M. Luttinger, {\em Lagrangian Tori in $\RR^4$}. J. Diff. Geom.  {52} (1999), 203--222.

 \bibitem{MS} D. McDuff and D. Salamon,
 `Introduction to symplectic topology.
Second edition.'  Oxford Mathematical Monographs. The Clarendon Press, Oxford University Press, New York, 1998. x+486


\bibitem{MMS} J. Morgan, T. Mrowka, and Z. Szabo, {\em Product formulas along $T\sp 3$ for Seiberg-Witten invariants.}  Math. Res. Lett. 4 (1997), no. 6, 915--929.

\bibitem{OS} B. Ozbaggi and A. Stipsitz, {\em Noncomplex smooth 4-manifolds with genus 2 Lefschetz fibrations}. Proc. Amer. Math. Soc. 128 (2000), no. 10, 3125--3128.

\bibitem{OzS} P. Ozsv\'{a}th and Z. Szab\'{o}, {\em On Park's exotic smooth four-manifolds,} `Geometry and topology of manifolds', 253--260, Fields Inst. Commun., { 47}, Amer. Math. Soc., Providence, RI, 2005.

\bibitem{OS2} P. Ozsv\'{a}th and Z. Szab\'{o}, {\em  The symplectic Thom conjecture.}  Ann. of Math. (2) 151 (2000), no. 1, 93--124.

\bibitem{dpark} B. Doug Park,{\em
Exotic smooth structures on $3\CP^2 \#n\bCP^2$.}
Proc. Amer. Math. Soc. 128 (2000), no. 10, 3057--3065.

\bibitem{park1} Jongil Park, {\em
Simply connected symplectic 4-manifolds with $b\sp +\sb 2=1$ and $c\sp 2\sb 1=2$. }
Invent. Math. 159 (2005), no. 3, 657--667.

\bibitem{park} Jongil Park, {\em The geography of symplectic 4-manifolds with arbitrary fundamental group.} Preprint (2006).

\bibitem{PSS} Jongil Park, A.  Stipsicz, and Z. Szabo, {\em
Exotic smooth structures on $\Bbb{CP}\sp 2\#5\overline{\Bbb{CP}\sp 2}$.}
Math. Res. Lett. 12 (2005), no. 5-6, 701--712.


\bibitem{perutz} T. Perutz, {\em A remark on K\"ahler forms on symmetric products of Riemann surfaces.}

\bibitem{Sal} D. Salamon, {\em Removable singularities and a vanishing theorem for Seiberg-Witten invariants,} Turkish J. Math {20} (1996), no. 1, 61-73.


\bibitem{smith} I. Smith, {\em Symplectic  geometry of Lefschetz fibrations}. Dissertation, Oxford 1998.

\bibitem{SS} A. Stipsicz and Z. Szab\'{o}, {\em An exotic smooth structure on $\CP^2\#6\overline{\CP}^2$,} Geom. Topol. {9} (2005), 813--832.

\bibitem{SS1} A. Stipsicz and Z. Szab\'{o}, {\em Small Exotic 4-manifolds with $b_2^+=3$,} Bull. London Math. Soc. {38} (2006), 501–506.

\bibitem{taubes1} C. Taubes, {\em Counting pseudo-holomorphic submanifolds in dimension  4.} J. Differential Geom. 44 (1996), no. 4, 818--893.

\bibitem{taubes2} C. Taubes, {\em Seiberg-Witten and Gromov invariants.} Geometry and physics (Aarhus, 1995), 591--601, Lecture Notes in Pure and Appl. Math., 184, Dekker, New York, 1997.

\bibitem{taubes3} C. Taubes, {\em The Seiberg-Witten invariants and symplectic forms.} Math. Res. Lett. 1 (1994), no. 6, 809--822.






\bibitem{usher} M. Usher, {\em Minimality and symplectic sums}. To appear in Internat. Math. Res. Not.  preprint (2006) math.SG/0606543

\end{thebibliography}
\end{document}